\documentclass[
final
]{dmtcs-episciences}

\usepackage[utf8]{inputenc}
\usepackage{subfigure}
\usepackage{amsthm, amsmath}

\usepackage[square]{natbib}

\newtheorem{theorem}{Theorem}[section]
\newtheorem{lemma}[theorem]{Lemma}
\newtheorem{claim}[theorem]{Claim}
\newtheorem{proposition}[theorem]{Proposition}

\newtheorem{corollary}[theorem]{Corollary}

\newtheorem{remark}[theorem]{Remark}

\author[Jovana Forcan and Jiayue Qi]{Jovana Forcan\affiliationmark{1,2}
  \and Jiayue Qi\affiliationmark{3}\thanks{ The research of the second author was funded by 
  the Austrian Science Fund (FWF): W1214-N15, project DK9.} }
  \title[ Maker–-Breaker domination number ]{Maker--Breaker domination number for Cartesian products of path graphs $P_2$ and $P_n$}
\affiliation{
  Department of Mathematics and Informatics, Faculty of Sciences, University of Novi Sad, Serbia\\
  Department of Computer Sciences and Systems, Faculty of Philosophy, University of East Sarajevo, Bosnia and Herzegovina\\
  Research Institute for Symbolic Computation, Johannes Kepler University, Linz, Austria}
\keywords{ Positional game,  Maker--Breaker domination game, Maker--Breaker domination number for grids, winning strategy on grids}
\begin{document}

\publicationdata{vol. 25:2 }{2023}{23}{10.46298/dmtcs.10465}{2022-12-13; 2022-12-13; 2023-06-22; 2023-10-19}{2023-10-22}
\maketitle
\begin{abstract}
 
We study the Maker--Breaker domination game played by Dominator and Staller on the vertex set of a given graph. Dominator wins when the 
vertices he has claimed form a dominating set of the graph. Staller wins if she makes it impossible for Dominator to win, or equivalently,
she is able to claim some vertex and all its neighbours. Maker--Breaker domination number $\gamma_{MB}(G)$ ($\gamma '_{MB}(G)$) of a graph 
$G$ is defined to be the minimum number of moves for Dominator to guarantee his winning when he plays first (second).
We investigate these two invariants for the Cartesian product of any two graphs. We obtain upper bounds for the Maker--Breaker domination number 
of the Cartesian product of two arbitrary graphs. Also, we give upper bounds for the Maker--Breaker domination number of the Cartesian product of
the complete graph with two vertices and an arbitrary graph. Most importantly, we prove that 
$\gamma'_{MB}(P_2\square P_n)=n$ for $n\geq 1$, $\gamma_{MB}(P_2\square P_n)$ 
equals~$n$, $n-1$, $n-2$, for $1\leq n\leq 4$, $5\leq n\leq 12$, and $n\geq 13$, respectively.
For the disjoint union of $P_2\square P_n$s, we show that 
$\gamma_{MB}'(\dot\cup_{i=1}^k(P_2\square P_n)_i)=k\cdot n$ ($n\geq 1$), and that
$\gamma_{MB}(\dot\cup_{i=1}^k(P_2\square P_n)_i)$ equals~$k\cdot n$, $k\cdot n-1$, $k\cdot n-2$ for 
$1\leq n\leq 4$, $5\leq n\leq 12$, and $n\geq 13$, respectively.
\end{abstract}
\section{Introduction}
\subsection{Background}
In this paper, we study the Maker--Breaker domination game which was first introduced in the 
literature by Duch\^{e}ne, Gledel, Parreau and Renault in~\citep{duchene2020maker}. 
This game combines two following research directions. In the original domination game introduced by Bre\v{s}ar, Klav\v{z}ar, 
and Rall in~\citep{BKR10}, two players, \textit{Dominator} and \textit{Staller},
alternately take a turn in claiming a vertex of the finite graph~$G$ which were not yet chosen in the course of the game --- note that 
the chosen vertex of either Dominator or Staller must enlarge the set of dominated vertices. 
Dominator has the goal to dominate the graph in as few moves as possible while Staller tries to prolong the game as much as possible. 
This indicates that both players play optimally in such a game. This is why later on, in our paper, we sometimes use 
subjunctive mood in the sentences, when we analyze the cases that are evidently not optimal for the player. 
However, in general, we keep the grammar in present tense, since many times it is not clear yet whether the focused moves are 
optimal for the player or not.

The Maker--Breaker games, introduced by Erd\H{o}s and Selfridge in~\citep{erdos1973combinatorial}, are played on a finite 
hypergraph $(X, \mathcal{F})$ with the vertex set~$X$ and a set $\mathcal{F} \subseteq 2^X$ of 
hyper-edges. The set~$X$ is called the \textit{board} of the game, while the set~$\mathcal{F}$ is called the \textit{family of winning sets}. 
Two players, \textit{Maker} and \textit{Breaker}, take turns in claiming previously unclaimed elements of $X$. Maker wins the game if, by
the end of the game, he has claimed all elements of some $F \in \mathcal{F}$.
Otherwise, Breaker wins. For a deeper and more comprehensive analysis of Maker--Breaker games, see the book of Beck~\citep{BeckBook08},
and the monograph of Hefetz, Krivelevich, Stojakovi\'{c} and Szab\'{o}~\citep{HKSSBook14}. 

The Maker--Breaker domination game (MBD for short) is played on a graph $G = (V,E)$ by two players, \textit{Dominator} and \textit{Staller}. 
The aim of Dominator is to build a dominating set of the graph, which is a set~$T$ such that every vertex not in~$T$ has a neighbour in~$T$. 
The aim of Staller is to claim a vertex from the graph $G$ and all its neighbours, so that Dominator cannot dominate this vertex, which leads 
to Dominator failing in dominating all vertices of~$G$.
As concluded in~\citep{duchene2020maker}, this game is equivalent to a Maker--Breaker game played on the vertex set of 
a given graph as the board, with winning sets being all closed neighbourhoods of vertices.

When it is not hard to determine the identity of the winner in some Maker--Breaker game, then a more interesting question to 
ask is how fast the player with the winning strategy can win.
Fast winning strategies for Maker in the Maker--Breaker games have received a lot of attention in the recent years 
(see for example~\citep{CFKL12, CFGHL15, HKSS09}). 
Specifically, for the Maker--Breaker domination game the smallest number of moves for Dominator is studied in~\citep{GIK19},
where Gledel, Ir\v{s}i\v{c}, and Klav\v{z}ar introduced the Maker--Breaker domination number $\gamma_{MB}(G)$ of a graph $G$,
as the minimum number of moves for Dominator to win in the game on $G$ where he is the first player. If Dominator is the second player,
then the corresponding invariant is denoted by $\gamma'_{MB}(G)$ in their paper.

\subsection{Preliminaries}\label{sec:pre}

Assume that the MBD game is in progress.
As in~\citep{GIK19}, we say that a game is the {\em D-game} 
if Dominator is the first to play, i.e. one \textit{round} consists of a move by Dominator followed by a move 
of Staller.
In the {\em S-game}, one round consists of a move by Staller followed by a move of Dominator. 
We denote by  $D_1, D_2,...$ (or $D'_1, D'_2,...$) the sequence of vertices chosen by
Dominator and by $S_1, S_2,...$ (or $S'_1, S'_2,...$) the sequence of vertices chosen by Staller, in a D-game (or in an S-game).
We say that the vertex~$v$ is \textit{isolated} by Staller if~$v$ and all its
neighbours are claimed by Staller. For a given graph~$G$, by $V(G)$ and $E(G)$ we denote its vertex set and edge set, respectively. 
If~$G$ is a graph and $S\subset V(G)$, then let~$G|S$ 
denote the graph~$G$ in which the vertices from~$S$ are declared to be already
dominated. 
As a preparation, we still need to introduce the {\em Continuation 
Principle} and the {\em No-Skip Lemma}.

\begin{remark}(The Continuation Principle, \citep[Remark~2.4]{GIK19})
\label{rem:cont-principle}
Let~$G$ be a graph with~$A,B\subset V(G)$. If $B\subset A$, then 
$\gamma_{MB}(G|A)\leq \gamma_{MB}(G|B)$ and 
$\gamma'_{MB}(G|A)\leq \gamma'_{MB}(G|B)$. 
\end{remark}

\begin{lemma}(No-Skip Lemma, \citep[Lemma~2.3]{GIK19})
\label{lem:no-skip}
In an optimal strategy of Dominator to achieve~$\gamma_{MB}(G)$
or $\gamma'_{MB}(G)$, it is never an advantage for him to skip
a move. Moreover, if Staller skips a move it can never disadvantage 
Dominator. 
\end{lemma}

\subsection{Main results}
In~\citep{GIK19}, the authors proposed finding the minimum number of moves for Dominator in the MBD game on 
the Cartesian product of two graphs. 
Motivated by the given problem, we give upper bounds for Maker--Breaker domination number for the 
Cartesian product of~$K_2$ and an arbitrary graph, and that 
of two general graphs as well,
in Section~\ref{sec:cartesian}. The corresponding results is Theorem~\ref{thm:bound1}. 
Most importantly, we focus on determining how fast can Dominator win on the graphs $P_2 \square P_n$, for $n\geq 1$.
From the results (Theorem~\ref{thm:grid_d} and Theorem~\ref{thm:grid_s}) on the grids, we get some results on the disjoint union 
of $P_2\square P_n$s as stated in Theorem~\ref{thm:union_d} and 
Theorem~\ref{thm:union_s}. We denote by $\dot\cup$ 
the disjoint union; and by $\dot\cup_{i=1}^k(G)_i$ we denote the disjoint union of~$k$ copies of graph~$G$.

The paper is organized as follows. 
In Section~\ref{sec:cartesian}, we prove Theorem~\ref{thm:bound1}. 
We prove Theorem~\ref{thm:grid_d} and Theorem~\ref{thm:grid_s} in Section~\ref{sec:p2_pn}. 
The corresponding proofs of Theorem~\ref{thm:union_d} and 
Theorem~\ref{thm:union_s} are given in Section~\ref{subsec:union}. We list the main results of this paper as the following.

\begin{theorem}\label{thm:bound1}
Let~$G$ and~$H$ be two arbitrary graphs on $n$ and $m$ vertices, respectively. Suppose that Dominator
has winning strategies on~$G$ and~$H$, for both D-games and S-games. Then
\[ \gamma _{MB}(G \square H) \leq \mathrm{min}\{\gamma _{MB}(G) + (m-1)\cdot \gamma' _{MB}(G), \gamma_{MB}(H) + (n-1)\cdot \gamma' _{MB}(H) \}.\]
Suppose that Dominator has a winning strategy for the S-games on both~$G$ and~$H$, then 
\[\gamma' _{MB}(G \square H) \leq \mathrm{min}\{m\cdot \gamma'_{MB}(G), n\cdot \gamma' _{MB}(H) \}.\]
\end{theorem}

\begin{theorem}\label{thm:grid_s}
For every positive integer~$n$, it holds that $\gamma' _{MB}(P_2 \square P_n) = n$. 
If Dominator skipped any moves during this game, he would lose the game. 
\end{theorem}

\begin{theorem}\label{thm:grid_d}
 Consider the D-games on $P_2\square P_n$. Then
 \begin{enumerate}
  \item If $1\leq n\leq 4$, then $\gamma_{MB}(P_2\square P_n)=n$.
  \item If $5\leq n\leq 12$, then $\gamma_{MB}(P_2\square P_n)=n-1$.
  \item If $n\geq 13$, then $\gamma_{MB}(P_2 \square P_n) = n-2$.
 \end{enumerate}
\end{theorem}

\begin{theorem}\label{thm:union_d}
For every two positive integers~$n$ and~$k$, it holds that 
\[\gamma_{MB}'(\dot\cup_{i=1}^k(P_2\square P_n)_i)=k\cdot n.\]
\end{theorem}

\begin{theorem}\label{thm:union_s}
Let~$k$ be any positive integer. Then
\begin{enumerate}
 \item If $1\leq n\leq 4$, then $\gamma_{MB}(\dot\cup_{i=1}^k(P_2\square P_n)_i) = k\cdot n$.
 \item If $5\leq n\leq 12$, then $\gamma_{MB}(\dot\cup_{i=1}^k(P_2\square P_n)_i) = k\cdot n-1$.
 \item If $n\geq 13$, then $\gamma_{MB}(\dot\cup_{i=1}^k(P_2\square P_n)_i) = k\cdot n-2$.
\end{enumerate}
\end{theorem}

\section{MBD games on the Cartesian product of two graphs}\label{sec:cartesian}
In this section, we prove the result on the Cartesian products of general graphs.

\begin{proof}[of Theorem~\ref{thm:bound1}]
Consider, first, the D-game on $G \square H$.
By $G^{(1)}, G^{(2)},...,G^{(m)}$ denote the $G$-layers of graph $G \square H$.
By $\mathsf{S}_D$ and $\mathsf{S}'_D$ denote Dominator's winning strategy on~$G$ in the D-game and in the 
S-game, respectively. 
Dominator will play his first move on one $G$-layer according to his winning strategy~$\mathsf{S}_D$.   
In every other round $i\geq 2$, he looks at the $(i-1)^{\mathrm{th}}$ move of Staller. 
If Staller in his $(i-1)^{\mathrm{th}}$ move claims a vertex from~$V(G^{(j)})$, let Dominator respond by claiming a vertex from the same set 
$V(G^{(j)})$ according to the corresponding winning strategy $\mathsf{S}_D$ or $\mathsf{S}'_D$ on graph~$G^{(j)}$ --- the choice
of $\mathsf{S}_D$ or $\mathsf{S}'_D$
depends on whether Dominator or Staller started 
the first move on this $G$-layer. Since Staller can be the first player on at 
most $m-1$ $G$-layers of graph $G \square H$, we know that Dominator can win within 
\[ \gamma _{MB}(G) + (m-1)\cdot\gamma' _{MB}(G)\] many moves. 

Since Dominator also has a winning strategy on~$H$ both as the first and as the second player, we obtain analogously that 
he can win within $\gamma_{MB}(H) + (n-1)\cdot \gamma' _{MB}(H)$ many moves. 
Therefore we get 
\[\gamma _{MB}(G \square H) \leq \mathrm{min}\{\gamma _{MB}(G) + (m-1)\cdot \gamma' _{MB}(G), \gamma_{MB}(H) + (n-1)\cdot \gamma' _{MB}(H) \}.\]

Now, we consider the S-game on $G\square H$. Staller starts the game, we let Dominator respond on the same $G$-layer, using 
the winning strategy he has on~$G$. Hence he can win within $m\cdot \gamma'_{MB}(G)$ many steps. 
Also, we can focus on the $n$ $H$-layers, in this Cartesian graph. With the analogous analysis, we know that 
Dominator can win within $n\cdot \gamma'_{MB}(H)$ many steps. Hence we get 
\[\gamma' _{MB}(G \square H) \leq \mathrm{min}\{m\cdot \gamma'_{MB}(G), n\cdot \gamma' _{MB}(H) \}.\]

Note that there can be a situation in any of the above cases, where Dominator cannot respond to the suggested strategy
since he has already dominated every vertex in the corresponding $G$-layer. 
In this case, either he skip the planned move or not does not influence the situation of the game.
The number we obtained above is in the situation where (we imagine that) he would not skip any move.
So, we know that, in reality, he needs even less moves to reach the winning status --- this is 
due to the ``No-Skip Lemma'' (Lemma~\ref{lem:no-skip}). Therefore, the above inequalities still hold,
under these situations.
\end{proof}

We have the following corollary from Theorem~\ref{thm:bound1}.
\begin{corollary}\label{cor:bound2}
Let $G$ be a graph on $n$ vertices. 
If Dominator has a wining strategy both as the first and as the second player in the game on $G$, then 
\[\gamma _{MB}(G \square K_2) \leq  \mathrm{min} \{\gamma _{MB}(G) + \gamma' _{MB}(G), n\}, \] and 
\[\gamma' _{MB}(G \square K_2) \leq  \mathrm{min} \{2\cdot \gamma' _{MB}(G), n\}.\]
\end{corollary}

\section{MBD game on $P_2\square P_n$}\label{sec:p2_pn}

In order to prove our main results on MBD games on grids, namely Theorem~\ref{thm:grid_s} and Theorem~\ref{thm:grid_d}, we need to introduce 
several ``MBD graphs'' and also two types of traps that Staller can make so as to win the game, as preparations.

\subsection{MBD graphs}

An {\em MBD graph} is a pair $(G,\mathcal{I})$, where $G=(V,E)$ is a graph 
and function 
\[\mathcal{I}\colon V\to \{s,d,n\}\times \{0,1\}\] assigns to each vertex a pair from $\{s,d,n\}\times \{0,1\}$, 
describing the current situation of the vertex.
Note that many MBD graphs defined in this section are standard partially dominated graphs.
However, some are not; also, we intend to have a consistent notation for all the graphs that show up in this section, 
and in precise mathematical language as well. Hence, we introduce this notion of ``MBD graphs''. 

To say it more intuitively,~$s$ means the vertex is already claimed by 
Staller,~$d$ means the vertex is already claimed by Dominator. And~$n$ refers to ``null'' which means that the vertex
is still free --- not yet claimed by any player in the current game. While~$1$ means that the vertex is already dominated,
i.e., it has a neighbouring vertex with the assigned value~$d$, and~$0$ means that the vertex is not yet dominated, i.e., 
none of its neighbours have been assigned to with the value~$d$ yet. 

An {\em MBD subgraph}~$S$ of a graph~$G$ is a subgraph of~$G$ which is an MBD graph, meaning 
that it has the value function $\mathcal{I}\colon V(S)\to \{s,d,n\}\times \{0,1\}$ attached to it.
In the sequel, we define several MBD subgraphs of the Cartesian graph $P_2\square P_m$. 
 They will be helpful in our later proofs. 
Let \[V_m = \{u_1,...,u_m, v_1,...,v_m\}\] and
\[E_m=\{\{u_i,u_{i+1}\} \mid i=1,2,...,m-1\} \cup \{\{v_i,v_{i+1}\}: i=1,2,...,m-1\}\]
\[\cup \{\{u_i,v_i\}: i=1,2,...,m\}.\]
It is easy to check that $P_2\square P_m=(V_m, E_m)$.

\begin{figure}[!h]
    \centering
    \subfigure[]{\includegraphics[width=0.4\textwidth]{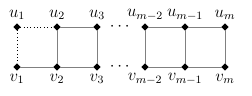}} 
    \hspace*{0.5cm}
    \subfigure[]{\includegraphics[width=0.4\textwidth]{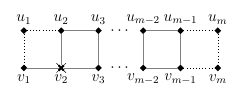}} 
    \newline
    \subfigure[]{\includegraphics[width=0.4\textwidth]{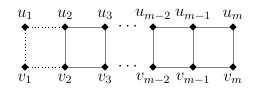}}
    \hspace*{0.5cm}
    \subfigure[]{\includegraphics[width=0.4\textwidth]{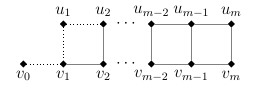}}
    \newline
    \subfigure[]{\includegraphics[width=0.4\textwidth]{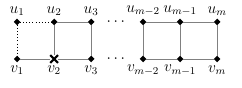}}
    \caption{Sub-figures: (a) $\mathfrak{X}_m$ (b) $\mathfrak{Y}_m$ (c) $\mathfrak{Z}_m$ (d) $\mathfrak{W}_m$ (e) $\mathfrak{R}_m$.
   Two incident edges being dotted indicates that the vertex is already dominated by Dominator; the cross indicates that 
   the vertex is claimed by Staller.}
    \label{fig:MBD_graphs}
\end{figure}

\begin{enumerate} 
\item[1.] By $\mathfrak{X}_m$ ($m\geq 1$) denote the MBD graph $(G,\mathcal{I})$ with $G = P_2 \square P_m$,~$\mathcal{I}$ assigns to vertex~$u_1$
with $(n,1)$ and to all other vertices with the value $(n,0)$. That is to say, vertex~$u_1$ is already dominated
in the ongoing game. Therefore, in the remaining game, Dominator does not need to consider dominating $u_1$. 
In the standard notation (see~\citep{GIK19}), $\mathfrak{X}_m$ is a partially dominated graph which is denoted by $G|\{u_1\}$.
(see Figure~\ref{fig:MBD_graphs}(a))
\item[2.] By $\mathfrak{Y}_m$ ($m\geq 3$) denote the MBD graph $(G,\mathcal{I})$ with $G= P_2 \square P_m$,~$\mathcal{I}$ assigns to vertex~$u_1$
value $(n,1)$, to vertex~$v_2$ value $(s,0)$, to vertices~$u_m$,~$v_m$ the value $(n,1)$ and to all other vertices the value $(n,0)$. 
When considering the D-game on~$\mathfrak{Y}_m$, we set $S_0 = v_2$ --- consider the MBD game played on this graph, keeping in mind that
Staller has claimed the vertex~$v_2$ and vertices~$u_m$ and~$v_m$ are already dominated (but not claimed by any players) before the game starts. 
(see Figure~\ref{fig:MBD_graphs}(b))
\item[3.] By $\mathfrak{Z}_m$ ($m\geq 1$) denote the MBD graph $(G,\mathcal{I})$ with $G= P_2 \square P_m$,~$\mathcal{I}$ assigns to vertices~$u_1$
and~$v_1$ the value $(n,1)$, to all other vertices value $(n,0)$. That is to say, vertices~$u_1$ and~$v_1$ are already dominated
in the ongoing game. Therefore, in the remaining game, Dominator does not need to consider dominating~$u_1$ or~$v_1$. 
In the standard notation (see~\citep{GIK19}),~$\mathfrak{Z}_m$ is a partially dominated graph which is denoted by $G|\{u_1, v_1\}$.
(see Figure~\ref{fig:MBD_graphs}(c))
\item[4.] By $\mathfrak{W}_m$ ($m\geq 1$) denote the MBD graph $(G,\mathcal{I})$, where $G=(V,E)$, and 
\[V=V_m\cup \{v_0\},\; E=E_m\cup \{\{v_0,v_1\}\};\]
$\mathcal{I}(v_0)=\mathcal{I}(u_1)=(n,1)$, and~$\mathcal{I}$ assigns to all other vertices value $(n,0)$. 
That is to say,~$v_0$ and~$u_1$ are already dominated 
before the game starts. 
In the standard notation (see~\citep{GIK19}), $\mathfrak{W}_m$ is a partially dominated graph which is denoted by $G|\{v_0,u_1\}$.
(see Figure~\ref{fig:MBD_graphs}(d))
\item[5.] By $\mathfrak{R}_m$ ($m\geq 0$) denote the MBD graph $(G,\mathcal{I})$ with $G=P_2\square P_m$ and 
\[\mathcal{I}(u_1)=(n,1),\; \mathcal{I}(v_2)=(s,0).\]
The function~$\mathcal{I}$ assigns to all other vertices value $(n,0)$. 
When considering the D-game on $\mathfrak{R}_m$, we set $S_0 = v_2$  --- consider the MBD game played on this graph, but Staller has claimed 
the vertex~$v_2$ and~$u_1$ is already dominated by Dominator before the game starts. Note that we use~$S_0$ to denote the 
``move'' of Staller before the game starts. Note that here we allow $m=0$, where~$\mathfrak{R}_0$ simply refers to the null graph; this 
allowance aims to make the later-on proof more smooth. When $m=1$, it just refers to the complete graph with two vertices.
(Figure~\ref{fig:MBD_graphs}(e)). 
\end{enumerate}

We define the {\em Maker--Breaker domination number of an MBD graph} --- abbreviated as ``domination number of an MBD graph'', in the later 
context of this paper
--- naturally as follows. Dominator intends to dominate all 
undominated vertices, namely those vertices with the second coordinate of~$\mathcal{I}(v)$ being~$0$; 
both players can only claim the unclaimed vertices, namely vertices~$v$ with the first coordinate of~$\mathcal{I}(v)$ being~$n$. We
use the same notations $\gamma_{MB}$ and $\gamma'_{MB}$ for the minimum number of moves for Dominator to guarantee his winning in 
D-game and S-game
on the MBD graph, respectively. Actually a normal graph can be viewed as a special case of an MBD graph, namely with all vertices assigned to with 
value $(n,0)$ by the function~$\mathcal{I}$.

Now we define two types of traps that Staller can create in the MBD game on $P_2 \square P_n$ for $n\geq 3$, so as to prevent the winning of Dominator. These two 
strategies of Staller will be used a lot later in the proofs.
\paragraph{Trap 1 - triangle trap.} We say that Staller created a \textit{triangle trap} if after her move Dominator is 
forced to claim the vertex $v_i$ (or $u_i$) --- so as to dominate $v_i$ (or $u_i$). This is
because all neighbouring vertices of $v_i$ (or $u_i$) are already claimed by Staller and she can 
isolate $v_i$ (or $u_i$) by claiming it in her next move,
if Dominator did not claim it. 
We say that Staller created a \textit{sequence of triangle traps} $v_i -- v_{j}$ (or $v_i -- u_j$), $i<j$, 
if Dominator is consecutively forced to claim vertices $v_i, u_{i+1}, v_{i+2}, u_{i+3},..., v_{j}$ (or $v_i, u_{i+1}, v_{i+2}, u_{i+3},..., u_{j}$)
because of the triangle traps one after another set up by Staller. 

\begin{figure}
    \centering
    \subfigure[]{\includegraphics[width=0.45\textwidth]{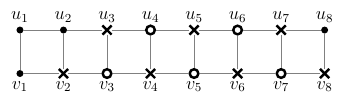}} 
    \hfill
    \subfigure[]{\includegraphics[width=0.45\textwidth]{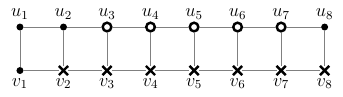}} 
    \caption{Illustrations of the two type of traps. (a) Sequence of triangle traps $v_3--v_7$; (b) Sequence of line traps $u_3--u_7$.}
    \label{fig:traps}
\end{figure}

The sequence of triangle traps $v_3--v_7$ is illustrated on Figure~\ref{fig:traps}(a). 
The first trap shows up when Staller has claimed $v_2$, $u_3$ and $v_4$, which forces Dominator to claim $v_3$; immediately after, Staller
claims $u_5$ to force Dominator to claim $u_4$, which creates the second trap. Notice that after Staller claims $v_8$, 
Dominator has to claim $v_7$; but Staller can then claim $u_8$, which will isolate $u_8$ and further leads 
to her winning of the game. This strategy will be applied a lot later on in our proofs.
The sequences of triangle traps $u_i--v_{j}$ and $u_i--u_j$ are defined analogously.

\paragraph{Trap 2 - line trap.} We say that Staller created a \textit{line trap} if after her move Dominator is forced to claim the vertex $v_i$ (or $u_i$) ---  
so as to dominate $u_i$ (or $v_i$). This is because all other neighbours of $u_i$ (or $v_i$) and $u_i$ (or $v_i$) itself are already claimed by Staller.
Staller can isolate $u_i$ (or $v_i$) by claiming $v_i$ ($u_i$) in her next move. 
We say that Staller creates a \textit{sequence of line traps} $v_i--v_j$ (or $u_i--u_j$), 
$i<j$, if Dominator is consecutively forced 
to claim the vertices $v_i, v_{i+1}, v_{i+2}, v_{i+3},..., v_{j}$ (or $u_i, u_{i+1}, u_{i+2}, u_{i+3},..., u_{j}$) 
because of the line traps one after another set up by Staller.

The sequence of line traps $u_3--u_7$ is illustrated on Figure~\ref{fig:traps}(b).
The first trap shows up when Staller has claimed $v_2$, $v_3$ and $v_4$, which forces Dominator to claim $u_3$ (in order to dominate $v_3$);
directly after, Staller claims $v_5$, which forces Dominator to claim $u_4$, which creates the second trap. Notice that after Dominator claims
$u_7$, Staller can then claim $u_8$, which makes it impossible for Dominator to dominate $v_8$ anymore, then win the game.
This strategy will be employed a lot later on in our proofs.

Next, let us briefly recall the {\em pairing strategy}. We follow the content provided 
in section~2 of~\citep{GIK19}. Let~$G$ be a graph, $k\geq 1$, and
$u_1,\ldots,u_k,v_1,\ldots,v_k$ pairwise different vertices of~$G$. 
Then we say that 
\[X=\{\{u_1,v_1\},\ldots,\{u_k,v_k\}\}\] is a 
{\em pairing dominating set} if 
\[\bigcup^k_{i=1}N[u_i,v_i]=V(G),\]
where~$N[u_i,v_i]$ means the set of common neighbouring vertices of~$u_i$
and~$v_i$. Note that here the set of ``neighbouring vertices'' of vertex~$v$ includes the vertex~$v$ itself.
When a graph has a pairing dominating set, then it is not hard to 
see that Dominator has a winning strategy in both S-game and D-game. He just 
need to claim the remaining vertex in the pair whenever there is one 
vertex left in a pair. We call this strategy the {\bf pairing strategy}.
Note that in this paper, the pairing strategy is adopted many times; most of the time it is 
adopted in order to prove the upper bound for the number of steps that Dominator needs in order to win. However, 
sometimes Dominator cannot conduct the game as suggested since the vertex that he ought 
to claim will not dominate any new vertices. In this case, we let Dominator skip that move
and continue the game. By the ``No-Skip Lemma'' (Lemma~\ref{lem:no-skip}), we know 
that Dominator --- in a real optimal game --- would need less steps to win. Hence the upper bound 
we tried to prove holds still. We only state this point once here, and omit it later on 
in the context where the pairing strategy are mentioned. 

With this we finished introducing the basic definitions which we will need later on in the proofs of our main results.

\subsection{Maker--Breaker domination numbers of MBD graphs}
In this section, we introduce some results on the Maker--Breaker domination numbers 
on MBD graphs, which then will lead to the completion of the proofs of Theorem~\ref{thm:grid_s} and Theorem~\ref{thm:grid_d}.
However, the proofs of Theorem~\ref{thm:grid_s} and Theorem~\ref{thm:grid_d} will be given in Section~\ref{sec:p2_pn}.
We first list the main results for this section, and then we give the proofs. 

\subsubsection{Domination numbers of MBD graphs: main results}
\begin{proposition}\label{prop:R_m}
  $\gamma_{MB}(\mathfrak{R}_m)=m$ for $m\geq 0$ --- note that~$\mathfrak{R}_0$ is simply the null graph; when $m\geq 2$, 
  Dominator would lose the game if he skipped any moves in the game.
 \end{proposition}
 
 \begin{proposition}\label{prop:y_m}
$\gamma _{MB}(\mathfrak{Y}_m) = m-1$, $m\geq 3$. 
\end{proposition}

\begin{proposition}\label{prop:z_m}
$\gamma' _{MB}(\mathfrak{Z}_m) = m-1$, $m\geq 1$. 
\end{proposition}

\begin{proposition}\label{prop:w_m}
 $\gamma' _{MB}(\mathfrak{W}_m) = m$ for $1\leq m\leq 3$, and $\gamma' _{MB}(\mathfrak{W}_m) = m-1$ for $m\geq 4$.
\end{proposition}

\begin{theorem}\label{thm:x_m}
  $\gamma_{MB}(\mathfrak{X}_1) = 1$, $\gamma_{MB}(\mathfrak{X}_m) = m-1$ for $2\leq m\leq 5$,
  and $\gamma _{MB}(\mathfrak{X}_m) = m-2$ for $m\geq 6$.
\end{theorem}

\subsubsection{Domination numbers of MBD graphs: proofs}

\begin{proof}[of Proposition~\ref{prop:R_m}]
First, the following five statements are not hard to verify:
 \begin{enumerate}
\item $\gamma_{MB}(\mathfrak{R}_0)=0$;
\item $\gamma_{MB}(\mathfrak{R}_1)=1$; 
\item $\gamma_{MB}(\mathfrak{R}_2)=2$;
\item Dominator would lose in the D-games on graph $\mathfrak{R}_2$ if he skipped
 any moves;
 \item On~$\mathfrak{R}_1$, Dominator would lose the S-game if he skipped any moves.
 \end{enumerate}
 
 Assume that $\gamma_{MB}(\mathfrak{R}_k)\leq k$ holds when $2\leq k\leq m-1$,
 then when $k=m$, let Dominator adopt the strategy for $\mathfrak{R}_{m-1}$
 on the (left) subgraph $\mathfrak{R}_{m-1}$, let him adopt the pairing strategy for the~$K_2$ (complete graph
 with two vertices) on the right. We then obtain that 
 \[\gamma_{MB}(\mathfrak{R}_m)\leq \gamma_{MB}(\mathfrak{R}_{m-1})+1\leq m.\]
 Hence by induction we obtain that
 $\gamma_{MB}(\mathfrak{R}_m)\leq m$ for any $m\geq 0$. 
 We only need to show that 
 $\gamma_{MB}(\mathfrak{R}_m)\geq m$ and Dominator would not win if he skipped any moves in the game.
 Assume that
 $\gamma_{MB}(\mathfrak{R}_k)\geq k$ and that Dominator would not win if he skipped any moves in D-game on graph $\mathfrak{R}_k$
 for any $k<m$. Now we consider the situation when $k=m$; notice that here we can assume 
 $m\geq 3$.
 
 In the case if Dominator had skipped the first move --- this situation is denoted by $D_1=\emptyset$ ---
 we propose the following strategy for Staller:
 Let $S_1=v_1$, which forces~$D_2$ to be~$u_1$. Then let Staller make line trap 
 from~$v_3$ to~$v_n$, which forces Dominator to respond on~$u_2$ until on~$u_{m-1}$. Then 
 let Staller claim~$u_m$, Dominator cannot dominate~$v_m$ anymore. Hence he would lose the game
 if Dominator skipped the first move.
 
 So now we consider all possibilities for~$D_1$. And we will propose correspondingly 
 Staller's strategy.
 \begin{itemize}
  \item Case 1: $D_1=u_i$, $i\geq 2$. Let $S_1=v_1$ which forces $D_2=u_1$. Then if $i\neq 2$,
  let Staller make line trap from~$v_3$ to~$v_i$, which forces Dominator
  to claim from~$u_2$ until~$u_{i-1}$. 
  If~$i=m$, then we already see that Dominator won with~$m$ moves.
  If~$i\notin\{m,m-1\}$, then let Staller claim~$v_{i+2}$, on the right 
  we obtain graph $\mathfrak{R}_{m-i}$. If $m-i=1$, just let Staller claim any remaining vertex.
  By induction we know that Dominator would not skip any moves and he needs~$m$ moves
  to win.
  
  \item Case 2: $D_1=v_i$, $i\geq 3$. Suppose $i>5$. Let $S_1=u_2$, which forces $D_2\in\{v_1,u_1\}$. 
  But no matter which choice~$D_1$ was assigned to, Staller will claim~$v_3$,~$u_3$ in the next two rounds so that either
  Dominator cannot dominate~$v_4$ or he cannot dominate~$u_4$. So he would lose the game, therefore this is not an optimal 
  choice for Dominator.
  Hence $i\in\{3,4\}$.
  \subitem Case 2.1: $i=3$. Staller will then claim~$u_1$ and~$u_3$, which forces Dominator
  to respond on~$v_1$ and~$u_2$, then Staller will claim~$u_5$, on the right we can use 
  the induction hypothesis. 
  \subitem Case 2.2: $i=4$. Let Staller claim~$u_2$, then Dominator needs to claim~$u_1$
  or~$v_1$.
  
  When $D_2=u_1$, then let Staller claim~$v_3$ which forces
  $D_3=v_1$. Then let $S_3=u_4$, which forces $D_4=u_3$. Then $S_4=u_6$, we can use induction
  hypothesis for the graph on the right. 
  
  When $D_2=v_1$, then let Staller claim~$u_3$ which forces $D_3=u_1$. Then let $S_3=u_4$ 
  which forces $D_4=v_3$. Then $S_4=u_6$, we can use the induction hypothesis for the graph on the right.
  
  \item Case 3: $D_1=u_1$, let $S_1=v_3$. In this case if $D_2=\emptyset$, Staller could then claim~$u_2$ 
  which would force Dominator to claim~$v_1$. Then let Staller claim~$u_3$, which gives Dominator a dilemma 
  of not being able to claim~$u_4$ and~$v_4$ at one step: he would fail. Hence Dominator would lose the game if skipped the move~$D_2$. For the 
  same reason, we know that $D_2\in \{v_1,u_2,u_3,u_4,v_4\}$. However if $D_2\in\{v_1,u_2\}$,
  then Staller could make line traps by claiming $v_4$ until $v_n$, then at the last step
  claiming~$u_n$, which would make Dominator lose the game. Hence $D_2\in \{u_3,u_4,v_4\}$.
  \subitem Case 3.1: $D_2=u_3$. Let $S_2=u_2$, forcing $D_3=v_1$. Then $S_3=v_5$, we can use the induction hypothesis
  for the graph on the right.
  \subitem Case 3.2: $D_2=u_4$. Let $S_2=u_2$, forcing $D_3=v_1$. Let $S_3=v_4$, forcing $D_4=u_3$.
  Then $S_4=v_6$, we can use the induction hypothesis for the graph on the right.
  \subitem Case 3.3: $D_2=v_4$. Let $S_2=u_2$, forcing $D_3=v_1$. Let $S_3=u_4$, forcing $D_4=u_3$.
  Then $S_4=u_6$, we can use induction for the graph on the right.
  
  \item Case 4: $D_1=v_1$, let $S_1=u_3$. In this case if $D_2=\emptyset$, Staller could then claim~$u_2$ which would force $D_3=u_1$.
  Then let Staller claim~$v_3$; we see that Dominator would not be able to dominate both~$u_3$ and~$v_3$ within only one step.
  So he would lose the game if he skipped~$D_2$. For the 
  same reason, we know that $D_2\in \{u_1,u_2,v_3,u_4,v_4\}$. However if $D_2\in\{u_1,u_2\}$,
  then Staller could make triangle traps by claiming~$v_4$ until $v_n/u_n$ (depending
  on whether~$m$ is odd or even), then at the last step
  claiming $u_m/v_m$, which would make Dominator lose the game. Hence $D_2\in \{v_3,u_4,v_4\}$.
  \subitem Case 4.1: $D_2=v_3$. Let $S_2=u_2$, forcing $D_3=u_1$. Then $S_3=u_5$, we can use the induction hypothesis
  for the graph on the right.
  \subitem Case 4.2: $D_2=u_4$. Let $S_2=u_2$, forcing $D_3=u_1$. Let $S_3=v_4$, 
  forcing $D_4=v_3$.
  Then $S_4=v_6$, we can use the induction hypothesis for the graph on the right.
  \subitem Case 4.3: $D_2=v_4$. Let $S_2=u_2$, forcing $D_3=u_1$. Let $S_3=u_4$, 
  forcing $D_4=v_3$.
  Then $S_4=u_6$, we can use the induction hypothesis for the graph on the right.
 \end{itemize}
 Notice that in some case it may happen that the graph is not big enough for the strategy
 of Staller, but it is not hard to verify that in those cases we simply have an~$\mathfrak{R}_0$ 
 or an~$\mathfrak{R}_1$ on the right, still we can use the induction hypothesis. 
 
 Hence notice that in all the above cases Dominator in total needs~$m$ moves to win and shall not skip 
 any moves during the game. 
 By induction, we conclude that for any $m\geq 2$, in the D-game on graph~$\mathfrak{R}_m$, 
 Dominator needs at least $m$
 steps to win and he shall not skip any moves in the game.
\end{proof}

\begin{remark}
Note that the D-game on graph~$\mathfrak{R}_m$ can be considered as the S-game on~$\mathfrak{X}_m$ with $S'_1 = v_2$. 
\end{remark}

\begin{proof}[of Proposition~\ref{prop:y_m}]
Vertex~$v_2$ is pre-claimed by Staller: we denote this by $S_0=v_2$.
The proof and case analysis can be done analogously to those of Proposition~\ref{prop:R_m}.
\end{proof}

\begin{proof}[of Proposition~\ref{prop:z_m}]
It is not hard to check that $\gamma' _{MB}(\mathfrak{Z}_m) = m-1$ holds for $m = 1,2,3$.
Assume that $\gamma'_{MB}(\mathfrak{Z}_k)\leq k-1$ for all $3 \leq k< m$. When $k=m$,  
we view graph~$\mathfrak{Z}_m$ as two parts --- $\mathfrak{Z}_{m-1}$, and~$u_m$,~$v_m$ together with the edge connecting them.
Of course these two parts are connected via two edges in~$\mathfrak{Z}_m$, namely $\{u_{m-1},u_m\}$ and $\{v_{m-1},v_m\}$ --- but we 
do not consider them for now, since the extra edges will just make it easier for Dominator to win. 
We propose the following strategy for Dominator. If Staller plays on the $\mathfrak{Z}_{m-1}$ part of the graph, 
let Dominator respond on $\mathfrak{Z}_{m-1}$ with his strategy on $\mathfrak{Z}_{m-1}$; if Staller plays on the other part, 
Dominator will use the pairing strategy, namely claim the other vertex in this part. In this way,
by the induction hypothesis, Dominator can win within 
\[\gamma'_{MB}(\mathfrak{Z}_{m-1})+1\leq (m-2)+1 = m-1\] steps. Hence,
\[\gamma'_{MB}(\mathfrak{Z}_m)\leq m-1,\; m\geq 1.\]

For the other direction of the proof, we need to show that $\gamma'_{MB}(\mathfrak{Z}_m)\geq m-1$ holds for $m\geq 4$.
We prove it by 
proposing a strategy for Staller. Let $S'_1=u_m$, then it has to be that $D'_1\in \{u_{m-1},v_{m-1},v_m\}$;
otherwise $S'_2=v_m$, then at least one of~$v_m$ and~$u_m$ can get isolated after~$S'_3$ (Staller's third step) --- 
after~$D'_2$ (Dominator's second step) at least one of $u_{m-1}$ and $v_{m-1}$ will still be free, so Staller can choose this free 
vertex at her third step. Now we make case distinctions for this three choices of~$D'_1$.
\begin{enumerate}
 \item $D'_1=v_{m-1}$. Let $S'_2=u_{m-1}$: this forces $D'_2=v_m$ (in order to dominate~$u_m$).
 Let $S'_3=u_{m-3}$. The remaining MBD graph in this game is a $\mathfrak{Y}_{m-2}$
 with vertex set 
\[V(\mathfrak{Y}_{m-2}) = \{u_1,u_2,...,u_{m-2}, v_1, v_2,...,v_{m-2}\}.\] It is Dominator's 
 turn now. Therefore he needs $\gamma_{MB}(\mathfrak{Y}_{m-2})=m-3$ (by Proposition~\ref{prop:y_m}) more moves to win. 
 So he needs in total $m-1$ steps to win in this case.
 \item $D'_1=u_{m-1}$. Let $S'_2=v_{m-1}$: this forces $D'_2=v_m$ (in order to dominate $v_m$). 
 Let $S'_3=v_{m-3}$. Remaining part of this case is the same as the last case. 
 Dominator needs in total $m-1$ steps to win.
 \item $D'_1=v_m$. Let $S'_2=u_{m-2}$: the remaining MBD graph in this game is a~$\mathfrak{Y}_{m-1}$.
 Dominator needs $\gamma_{MB}(\mathfrak{Y}_{m-1})=m-2$ (by Proposition~\ref{prop:y_m}) many moves to win on this part and he needs
 in total $(m-2)+1= m-1$ many moves to win in this case. 
\end{enumerate}
Hence, Staller has a strategy such that Dominator needs at least $m-1$ steps in order to win in an 
S-game on the MBD graph $\mathfrak{Z}_m$. We obtain that $\gamma'_{MB}(\mathfrak{Z}_m)\geq m-1$, $m\geq 4$.
\end{proof}

\begin{proof}[of Proposition~\ref{prop:w_m}]
One can check that 
\[\gamma' _{MB}(\mathfrak{W}_m) = m,\; 1\leq m\leq 3.\] Let $m\geq 4$.
Since $\mathfrak{W}_m$ has one more undominated vertex than $\mathfrak{Z}_m$ (namely $v_1$), 
Dominator needs to play at least as many moves on $\mathfrak{W}_m$ as on $\mathfrak{Z}_m$.
Therefore we have \[\gamma'_{MB}(\mathfrak{W}_m) \geq \gamma'_{MB}(\mathfrak{Z}_m)=m-1,\] which shows the lower bound. 
For the upper bound, we start our consideration from $m=4$. 
When $m=4$, we propose strategies for Dominator accordingly to the steps of Staller.
\begin{itemize}
\item Case 1: $S'_1 = v_2$.  Let $D'_1 = u_3$. 
\subitem If $S'_2 = v_1$, then let $D'_2 = v_3$. Dominator only needs one step to dominate~$v_1$ before his winning; so as to do this, he only needs to 
claim~$v_0$ or~$u_1$ in the next move. So he needs in total~$3$ moves to win.
\subitem If $S'_2\neq v_1$, then let $D'_2=v_1$. Then Dominator just needs one more step to dominate~$v_4$. In total he needs~$3$ steps to win.

\item Case 2: $S'_1 \neq v_2$. 
\subitem If $S'_1=u_4$ (or $v_4$). Let $D'_1 = u_3$. 
If $S'_2 = v_3$, then let $D'_2 = v_4$ (or $u_4$) and let $D'_3\in \{v_1,v_2\}$. 
If $S'_2 = v_4$ (or $u_4$), 
let $D'_2 = v_3$ and $D'_3\in \{v_1,v_2\}$. So Dominator needs in total~$3$ steps to win.
\subitem  If $S'_1 \notin \{u_4,v_4\}$. 
Let $D'_1 = v_2$. Dominator needs at most two more moves to dominate the remaining vertices. 

\end{itemize} 
So far, we obtain that 
$\gamma'_{MB}(\mathfrak{W}_4)\leq 3$. Suppose that \[\gamma'_{MB}(\mathfrak{W}_{m-1})\leq m-2\; \text{for}\; m\geq 5.\]
Then consider the MBD graph $\mathfrak{W}_m$ ($m\geq 5$) as two parts (of course plus the edges connecting the two parts): 
one is the rightmost~$K_2$, namely~$u_m$,~$v_m$ 
and the edge connecting them; the other is on the left an MBD subgraph $\mathfrak{W}_{m-1}$. We let Dominator react 
on the left $\mathfrak{W}_{m-1}$ according to his strategy on $\mathfrak{W}_{m-1}$ whenever Staller claims a vertex on this subgraph; otherwise
we let Dominator claim the remaining vertex among~$u_m$ and~$v_m$. By induction, 
Dominator needs no more than $(m-2)+1=m-1$ moves to win. Hence $\gamma'_{MB}(\mathfrak{W}_m)\leq m-1$. 
\end{proof}

\begin{figure}
\centering
\includegraphics[width=0.74\linewidth]{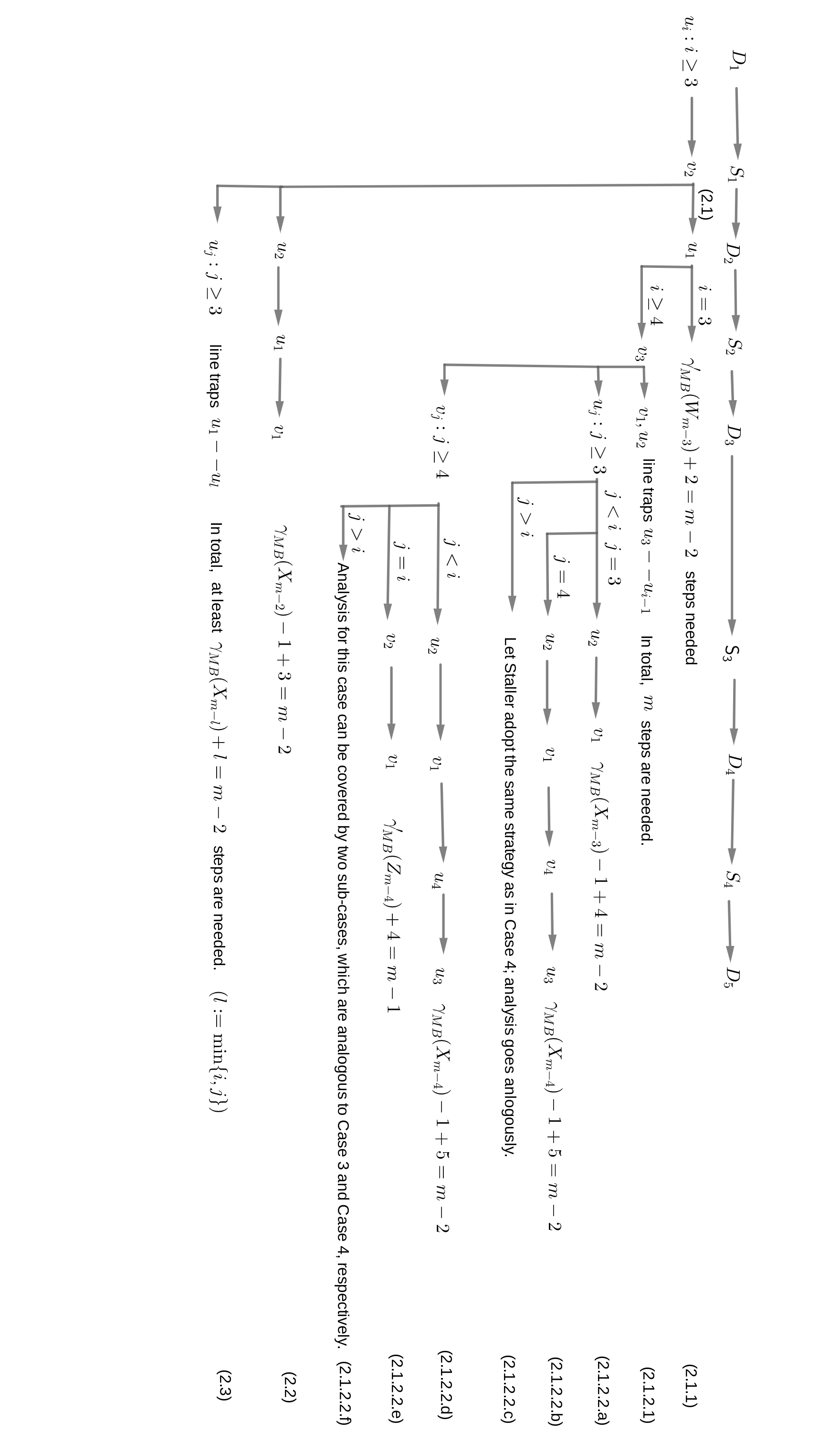}
\caption{This is a strategy tree showing the strategies for Staller under the situation when 
$D_1=u_i(i\geq 3)$, $S_1=v_2$, and $D_3=u_i(1\leq i\leq m)$, in an MBD game on~$\mathfrak{X}_m$.
At the end of each branch, we explain at least how many steps Dominator would need in order to win the game. It covers 
the cases (2.1), (2.2) and (2.3) in the proof of Theorem~\ref{thm:x_m}.}
\thispagestyle{empty}
\label{fig:strategy_tree}
\end{figure}

Theorem~\ref{thm:x_m} is the essence for proving Theorem~\ref{thm:grid_d}. The case analyses in the proof can be lengthy.
We provide a strategy tree (see Figure~\ref{fig:strategy_tree}) showing the idea of the case studies
of Theorem~\ref{thm:x_m}, which, although not covering all cases, 
but can at least express the idea behind. We leave it as an exercise for readers, to complete the whole strategy tree taking into consideration 
of all cases.

\begin{proof}[of Theorem~\ref{thm:x_m}]
For $m\in \{1,2,3\}$ the situation is not hard to directly see. 
For $m=4$ and $m=5$ simple cases analysis gives the result. 
Let $m\geq 6$. First, we consider the D-game on $\mathfrak{X}_6$.  Let $D_1 = v_2$, then we see an MBD subgraph $\mathfrak{W}_4$ on the rightmost which already 
contains all the undominated vertices. 
By Proposition~\ref{prop:w_m}, we know that $\gamma'_{MB}(\mathfrak{W}_4) = 3$. Hence, $\gamma_{MB}(\mathfrak{X}_6) =  4$. 

Suppose that $\gamma_{MB}(\mathfrak{X}_k)\leq k-2$ for $k<m$. When $k=m$, consider the graph as two parts (plus edges connecting the two parts):
one is the $K_2$ on the rightmost, the other is the MBD subgraph $\mathfrak{X}_{m-1}$ on the left. Let Dominator 
adopt the strategy for $\mathfrak{X}_{m-1}$ on the left subgraph whenever Staller plays on that subgraph and let him claim the remaining vertex 
among $u_m$ and $v_m$ when Staller claims one of them. In this way, by the induction hypothesis we know that 
Dominator needs at most $(m-3)+1=m-2$ steps to win. Hence $\gamma_{MB}(\mathfrak{X}_m)\leq m-2$ for $m\geq 6$. 
 
We prove the lower bound by induction and proposing strategies for Staller. 
Suppose that 
\[\gamma_{MB}(\mathfrak{X}_k)\geq k-2,\; \text{ for }\; 4\leq k< m,\; m\geq 6.\]
Now consider the case when $k=m$, $m\geq 6$. We want to show that \[\gamma_{MB}(\mathfrak{X}_m)\geq m-2,\; m\geq 6.\]
According to the first step of Dominator, there are three cases: $D_1\in \{u_1,v_1,u_2,v_2\}$, $D_1=u_i$ ($i\geq 3$), or $D_1=v_i$ ($i\geq 3$).

\begin{itemize}
 \item{\bf Case 1:} $D_1\in \{u_1,v_1,u_2,v_2\}$. 
 \subitem If $D_1\in \{u_1,v_1\}$, then we see that the remaining game is the S-game on the MBD subgraph $\mathfrak{W}_{m-1}$. By Proposition~\ref{prop:w_m}, 
 $\gamma'_{MB}(\mathfrak{W}_m)=m-1$ when $m\geq 4$. Hence Dominator needs $m-2$ more steps to win; he needs in total $m-1$ steps to win.
 \subitem If $D_1=v_2$, then we see that on the right is an MBD subgraph $\mathfrak{W}_{m-2}$. By Proposition~\ref{prop:w_m}, 
 Dominator needs in total $(m-3)+1=m-2$ moves to win. 
 \subitem If $D_1=u_2$, then on the left is the MBD subgraph $\mathfrak{W}_{m-2}$, while on the right Dominator 
 needs at most one more step, so as to dominate $v_1$. Hence Dominator needs at least $(m-3)+1=m-2$ steps to win. 
 
 \item{\bf Case 2:} $D_1 = u_i$, $i\geq 3$. Let $S_1=v_2$, then according to the second move of Dominator, we carry out case distinctions.
 
 \paragraph{(2.1)} $D_2=u_1$. Then we further consider two sub-cases: $i=3$ or $i\geq 4$. 
 
 \paragraph{(2.1.1)} If $D_1=u_3$, then we see that the remaining game is the S-game on the MBD subgraph $\mathfrak{W}_{m-3}$ on the left of $\mathfrak{X}_m$; where on the right
          Dominator needs at most one more step, so as to dominate $v_2$. Hence he needs in total at least $(m-4)+2=m-2$ steps to win. 
   \paragraph{(2.1.2)} If $D_1=u_i$, $i\geq 4$, we let $S_2=v_3$. We need to further distinct cases according to the choice of $D_3$. 
   
    \paragraph{(2.1.2.1)} $D_3\in\{v_1,u_2\}$. Let Staller take the strategy of creating a sequence of line traps $u_3--u_{i-1}$ which ends with
    $S_{i-1}=v_i$ and forces $D_i= u_{i-1}$. Then, if $m-i\geq 2$, let $S_i=v_{i+2}$. We see that the rightmost is the MBD graph 
    $\mathfrak{R}_{m-i}$. Dominator wins after the D-game on the rightmost $\mathfrak{R}_{m-i}$. 
    By Proposition~\ref{prop:R_m}, $\gamma_{MB}(\mathfrak{R}_{m-i})=m-i$. Hence Dominator needs in total $m$ moves to win. If 
    $m-i=1$, then let $S_i\in\{u_m,v_m\}$, Dominator needs only one more step to win. If $m-i=0$, the game finished already.
    In either case, Dominator wins with $m$ steps. 
   \paragraph{(2.1.2.2)} $D_3=u_j$, $j\geq 3$ or $D_3=v_j$, $j\geq 4$. Suppose that $\min\{i,j\}\notin\{3,4\}$. Let $S_3=u_2$, which forces
    $D_4=v_1$. Then let $S_4=u_3$, Dominator cannot dominate both $u_4$ and $v_4$, hence he would not be able to win. Therefore, 
     $\min\{i,j\}\in\{3,4\}$.
     \paragraph{(2.1.2.2.a)} When $D_3= u_j$, where $j\geq 3$ and $j<i$, we know that $j\in\{3,4\}$. 
      When $j=3$, we have $D_3=u_3$. Recall that $D_1=u_i$, $i\geq 4$, $S_1=v_2$, $D_2=u_1$, $S_2=v_3$, $D_3=u_3$.
      Now let $S_3=u_2$, which forces $D_4=v_1$. We see on the rightmost the MBD graph $\mathfrak{X}_{m-3}$, but, with 
      $u_i$ being claimed by Dominator already. Now it is Staller's turn. By the induction hypothesis, we know that 
      Dominator needs at least $\gamma_{MB}(\mathfrak{X}_{m-3})-1\geq m-6$ more steps to win the game. In total, he needs 
      $(m-6)+4=m-2$ steps.
      
       \paragraph{(2.1.2.2.b)} When $j=4$, we have $D_3=u_4$. Let $S_3=u_2$, which forces $D_4=v_1$. Then let $S_4=v_4$, which forces 
      $D_5=u_3$. Similarly as in the last case, the remaining graph is the MBD graph $\mathfrak{X}_{m-4}$ on the right, but
      with $u_i$ being claimed by Dominator already. Now it is Staller's turn. By the induction hypothesis, we know that 
      Dominator needs at least $\gamma_{MB}(\mathfrak{X}_{m-4})-1\geq m-7$ more steps to win the game. In total, he needs 
      $(m-7)+5=m-2$ steps.
      
      \paragraph{(2.1.2.2.c)} When $D_3= u_j$, where $j>i$, we know that $i\in\{3,4\}$. Since $i\geq 4$, we know that $D_1=u_i=u_4$.
      Let Staller take the same strategy as in the last case, we get that Dominator needs at least $m-2$ steps to win.
      
      \paragraph{(2.1.2.2.d)} When $D_3=v_j$, where $j<i$, since $S_2=v_3$, we know that $j=4$, i.e., $D_3=v_4$. Let $S_3=u_2$,
      which forces $D_4=v_1$. Then let $S_4=u_4$, which forces $D_5=u_3$. The remaining graph on the right is the 
      MBD graph $\mathfrak{X}_{m-4}$, but with $u_i$ claimed by Dominator already, and now it is Staller's turn. By the 
      induction hypothesis, we know that 
      Dominator needs at least $\gamma_{MB}(\mathfrak{X}_{m-4})-1\geq m-7$ more steps to win the game. In total, he needs 
      $(m-7)+5=m-2$ steps.
      
      \paragraph{(2.1.2.2.e)} When $D_3=v_j$, where $j=i$, since $S_2=v_3$, we know that $j=i=4$, that is, $D_1=u_4$ and $D_3=v_4$.
      Let $S_3=v_2$, which forces $D_4=v_1$. The remaining graph on the right is the MBD graph $\mathfrak{Z}_{m-4}$, and now it is 
      Staller's turn. By Proposition~\ref{prop:z_m}, Dominator needs in $\gamma'_{MB}(\mathfrak{Z}_{m-4})=m-5$ more steps to win. 
      In total, he needs $(m-5)+4=m-1$ steps to win.
      
      \paragraph{(2.1.2.2.f)} When $D_3=v_j$, where $j>i$, then $i\in \{3,4\}$, that is, $D_1\in\{u_3,u_4\}$. The analyses are
      analogous to items (a) and (b). Dominator needs in total at least $m-2$ steps to win.

      \paragraph{(2.2)} $D_2=u_2$. 
    Recall that we are currently under the case where $D_1=u_i$, $i\geq 3$, $S_1=v_2$. Let $S_2=u_1$, which forces $D_3=v_1$. 
   The remaining graph is the MBD graph $\mathfrak{X}_{m-2}$ on the right, but with $u_i$ claimed already by Dominator, 
   and now it is Staller's turn. By the induction hypothesis,
   Dominator needs at least $m-5$ more steps to win. In total, he needs at least $(m-5)+3=m-2$ steps to win.
   
   \paragraph{(2.3)} $D_2=u_j$, where $j\geq 3$. 
   Recall that we are currently under the case where $D_1=u_i$, $i\geq 3$, 
   $S_1=v_2$. Now let $S_2=v_1$, which forces $D_3=u_1$. Let $l:=\min\{i,j\}$ and $h:=\max\{i,j\}$. Let Staller create a 
   sequence of line traps $u_1--u_{l-1}$. The remaining graph is the MBD graph $\mathfrak{X}_{m-l}$ on the right, but with $u_h$ claimed already
   by Dominator, and now it is Staller's turn. By induction, Dominator needs at least $m-l-2$ more steps to win. 
   In total, he needs at least $(m-l-2)+l=m-2$ steps 
   to win.

\paragraph{(2.4)} $D_2=v_1$.

\paragraph{(2.4.1)} $D_1=u_3$. Recall that $S_1=v_2$, $D_2=v_1$. The remaining graph is the MBD graph $\mathfrak{W}_{m-3}$ on the right, 
 and now it is Staller's turn. By Proposition~\ref{prop:w_m}, we know that Dominator needs $\gamma'_{MB}(\mathfrak{W}_{m-3})\geq m-4$
 more steps to win. In total, he needs at least $(m-4)+2=m-2$ steps to win.
 
 \paragraph{(2.4.2)} $D_1=u_i$, where $i\geq 4$. Let $S_2=u_3$. We need to further distinct cases according to the choice of $D_3$. 
  
    \paragraph{(2.4.2.1)} $D_3\in\{u_1,u_2\}$. 
    
    \paragraph{(2.4.2.1.a)} When $i$ is even, let Staller create a sequence of 
    triangle traps $v_3--v_{i-1}$, which ends with $S_{i-1}=v_i$, $D_i=v_{i-1}$. Then, let $S_i=v_{i+2}$ if 
    $m-i\geq 2$. The remaining graph is the MBD graph $\mathfrak{R}_{m-i}$ and now it is Dominator's turn. 
    By Proposition~\ref{prop:R_m}, we know that Dominator needs $\gamma_{MB}(\mathfrak{R}_{m-i})=m-i$ more steps to win.
    In total, he needs $(m-i)+i=m$ steps to win.
    \paragraph{(2.4.2.1.b)} When $i$ is odd, let Staller create a sequence of triangle traps $v_3--v_{i-2}$, which ends with 
    $S_{i-2}=v_{i-1}$, $D_{i-1}=v_{i-2}$. The remaining graph is the MBD graph $\mathfrak{W}_{m-i}$ on the right and now 
    it is Staller's turn. By Proposition~\ref{prop:w_m}, Dominator needs at least $\gamma'_{MB}(\mathfrak{W}_{m-i})\geq m-i-1$
    more steps to win. In total, he needs at least $(m-i-1)+i-1=m-2$ steps to win.

    \paragraph{(2.4.2.2)} $D_3=u_j$, $j\geq 4$ or $D_3=v_j$, $j\geq 3$. A similar reasoning as in case {\bf (2.1.2.2)}
    shows that $\min\{i,j\}\in\{3,4\}$ must hold.
    
     \paragraph{(2.4.2.2.a)} When $D_3= u_j$, where $j<i$; since $S_2=u_3$, we know that $j=4$, that is,
   $D_3=u_4$. Let $S_3=u_2$, which forces $D_4=u_1$; then let $S_4=v_4$, which forces $D_5=v_3$. 
   The remaining graph is the MBD graph $\mathfrak{X}_{m-4}$, but with $u_j$ already claimed by Dominator, and now it is Staller's turn.
   By the induction hypothesis, Dominator needs at least $m-6-1=m-7$ more steps to win. In total, he needs 
   $(m-7)+5=m-2$ steps to win.

      \paragraph{(2.4.2.2.b)} When $D_3= u_j$, where $j>i$; since $S_2=u_3$, we know that $i=4$, i.e., $D_1=u_4$. 
      With the same strategy as in item (a), we obtain that Dominator needs at least $m-2$ steps to win.
      
      \paragraph{(2.4.2.2.c)} When $D_3=v_j$, where $j>i$; since $S_2=u_3$, we know that $i=4$, i.e., $D_1=u_4$. 
      With the same strategy as in item (a), we obtain that Dominator needs at least $m-2$ steps to win. 
      
      \paragraph{(2.4.2.2.d)} When $D_3=v_j$, where $j=i$; since $S_2=u_3$, we know that $i=j=4$, i.e., $D_1=u_4$, $D_3=v_4$.
      Let $S_3=u_2$, which forces $D_4=u_1$. The remaining graph is the MBD graph $\mathfrak{Z}_{m-4}$ on the right and it is 
      Staller's turn.
      By Proposition~\ref{prop:z_m}, Dominator needs at least $\gamma'_{MB}(\mathfrak{Z}_{m-4})=m-5$ more steps to win. 
      
      \paragraph{(2.4.2.2.e)} When $D_3=v_j$, where $j<i$, then $j\in\{3,4\}$. When $j=3$, $D_3=v_3$. Let $S_3=u_2$,
      which forces $D_4=u_1$. Consider the remaining graph as the MBD graph $\mathfrak{X}_{m-3}$, but with $u_i$ already
      claimed by Dominator, and now it is Staller's turn. By the induction hypothesis, Dominator needs at least $k-6$ more steps to win. In total, he needs
      at least $k-2$ steps to win. When $j=4$, $D_3=v_4$. Let $S_3=u_2$, which forces $D_4=u_1$. Then let 
      $S_4=u_4$, which forces $D_5=v_3$. We get the MBD graph $\mathfrak{X}_{m-4}$ with $u_i$ already claimed by Dominator on 
      the right. By induction, Dominator needs at least $k-7$ steps to win. In total, he needs at least 
      $k-2$ steps to win.

  \paragraph{(2.5)} $D_2=v_j$, where $j>i$.
   
   Recall that we are currently under the case where $D_1=u_i$, $i\geq 3$, $S_1=v_2$. Let $S_2=v_1$, which forces
   $D_3=u_1$. Then let Staller create a sequence of line traps $u_2--u_{i-1}$. Then the remaining graph is the MBD 
   graph $\mathfrak{X}_{m-i}$ on the right, but with $v_j$ already claimed by Dominator, and now it is Staller's turn. 
   By induction, Dominator 
   needs at least $m-i-3$ more steps to win. In total, he needs $(m-i-3)+i+1=m-2$ steps to win.

    \paragraph{(2.6)} $D_2=v_j$, where $j=i\geq 3$.
   
   Recall that $D_1=u_i$, $i\geq 3$, $S_1=v_2$, $D_2=v_i$. Let $S_2=v_1$, which forces $D_3=u_1$. Then let Staller
   create a sequence of line traps $u_2--u_{i-2}$. If $i=m$, then Dominator already won with~$m$ steps; otherwise,
   the remaining graph is the MBD graph $\mathfrak{Z}_{m-i}$ on the right. By Proposition~\ref{prop:z_m}, Dominator needs 
   $\gamma'_{MB}(\mathfrak{Z}_{m-i})=m-i-1$ more steps to win. In total, he needs $(m-i-1)+i=m-1$ steps to win. In either case,
   he needs more than $m-2$ steps to win. 
   
    \paragraph{(2.7)} $D_2=v_j$, where $i>j\geq 2$ and $j$ is even. 
   
   Recall that $D_1=u_i$, $i\geq 3$, $S_1=v_2$, so actually $j\geq 4$. Let $S_2=u_2$, then we claim that $D_3\in \{u_1,v_1\}$. 
   Suppose not. After $D_3$ either $v_3$ or $u_3$ should be free. If $v_3$ is free after $D_3$, let 
   $S_3=v_1$ and $S_4\in\{u_1,v_3\}$. We see that Dominator cannot dominate both $u_1$ and $v_2$, hence he would not win. 
   If $u_3$ is free after $D_3$, let $S_3=u_1$ and $S_4\in\{v_1,u_3\}$. 
   We see that Dominator cannot dominate both $v_1$ and $u_2$, hence he would not win. 
   
   \paragraph{(2.7.a)} $D_3=u_1$. Recall that $D_1=u_i$, $i\geq 3$, $S_1=v_2$, $D_2=v_j$, $i>j\geq 4$ is even, $S_2=u_2$.
   Let $S_3=v_3$, which forces $D_4=v_1$. Then let Staller create a sequence of triangle traps $u_3--u_{j-1}$. 
   The remaining graph is the MBD graph $\mathfrak{X}_{m-j}$ on the right, but with $u_i$ claimed already by Dominator, and 
   now it is Staller's turn. By induction, Dominator needs at least $m-j-2$ more steps to win.
   In total, he needs $(m-j-2)+j=m-2$ steps to win. 
   
   \paragraph{(2.7.b)} $D_3=v_1$. Recall that $D_1=u_i$, $i\geq 3$, $S_1=v_2$, $D_2=v_j$, $i>j\geq 4$ is even, $S_2=u_2$.
   Let $S_3=u_3$, which forces $D_4=u_1$. Actually $j=4$ must hold, otherwise in the next move Staller can either isolate
   $u_4$ or $v_4$. Let $S_4=u_4$, which forces $D_5=v_3$. The remaining graph is the MBD graph $\mathfrak{X}_{m-4}$
   on the right, but with $u_i$ already claimed by Dominator, and now it is Staller's turn. By induction, Dominator needs at least 
   $m-7$ more steps to win. In total, he needs at least $(m-7)+5=m-2$ steps to win.
   
    \paragraph{(2.8)} $D_2=v_j$, where $i>j\geq 2$ and $j$ is odd. Let $S_2=u_1$, which forces $D_3=v_1$.
   Then let Staller create a sequence of triangle traps $u_2--u_{j-1}$. The remaining graph is the MBD graph $\mathfrak{X}_{m-j}$
   on the right, but with $u_i$ already claimed by Dominator, and now it is Staller's turn. By induction, Dominator needs at least 
   $m-j-3$ more steps to win. In total, he needs at least $(m-j-3)+j+1=m-2$ steps to win.
   
   \item {\bf Case 3:} $D_1=v_i$, $i\geq 3$.  
   Let $S_1=v_2$, then according to the second move of Dominator, we carry out case distinctions.
   
    \paragraph{(3.1)} $D_2=u_1$.
    
    \paragraph{(3.1.1)} If $i=3$, i.e., $D_1=v_3$, then the remaining graph is the MBD graph $\mathfrak{W}_{m-3}$
    on the right, and it is Staller's turn now. By Proposition~\ref{prop:w_m}, Dominator needs
    at least $(m-3)-1=m-4$ more steps to win. In total, he needs $(m-4)+2=m-2$ steps to win. 
    
    \paragraph{(3.1.2)} If $i\geq 3$, let $S_2=v_3$. We carry out case distinctions according to the third 
    move of Dominator in the sequel. 
    
    \paragraph{(3.1.2.1)} $D_3=v_1$. If $i=4$, i.e., $D_1=v_4$, on the right we see the MBD graph $\mathfrak{W}_{m-4}$,
    and it is Staller's turn now. Dominator needs at least $(m-4)-1=m-5$ more steps to dominate 
    those undominated vertices in this graph $\mathfrak{W}_{m-4}$. Hence he needs in total at least $(m-5)+3=m-2$
    steps to win. Otherwise, if $i>4$, let $S_3=v_4$ starting a sequence of line traps $u_3--u_{i-2}$.
    The remaining graph is the MBD graph $\mathfrak{W}_{m-i}$ on the right, and it is Staller's turn now.
    By Proposition~\ref{prop:w_m}, Dominator needs at least $(m-i)-1$ more steps to win. In total, 
    he needs at least $(m-i-1)+(i-1)=m-2$ steps to win.  
    
    \paragraph{(3.1.2.2)} $D_3=v_j$, $j\geq 4$ or $D_3=u_j$, $j\geq 3$. A similar reasoning as in case {\bf (2.1.2.2)}
    shows that $\min\{i,j\}\in\{3,4\}$ must hold. 
    
    \paragraph{(3.1.2.2.1)} $D_3=v_j$, $j\geq 4$. Denote by $l:=\min\{i,j\}$, $h:=\max\{i,j\}$, then we know that $l=4$ since $v_3$ is claimed 
    by the second move of Staller. Let $S_3=u_2$, which forces $D_4=v_1$. Then let $S_4=u_4$, which forces $D_5=u_3$ 
    by making a triangle trap.
    The remaining graph is the MBD graph $\mathfrak{X}_{m-4}$ on the right, but with $u_h$ claimed by Dominator, and it is Staller's turn now.
    By induction, Dominator needs at least $(m-4)-2=m-6$ more steps to win. In total, he needs
    at least $(m-6)+4=m-2$ steps to win. 
    
    \paragraph{(3.1.2.2.2)} $D_3=u_j$, $j>i$. Since $v_3$ is already claimed, we know that $i=4$, that is, $D_1=v_4$. 
    The reasoning for this case is analogous to the last case (case {\bf (3.1.2.2.1)}).
    
    \paragraph{(3.1.2.2.3)} $D_3=u_j$, $j=i$. Since $v_3$ is already claimed, we know that $i=j=4$. The MBD graph on the 
    right is $\mathfrak{Z}_{m-4}$, and it is Staller's turn now. By Proposition~\ref{prop:z_m}, Dominator needs at least 
    $(m-4)-1=m-5$ more steps to dominate those undominated vertices in $\mathfrak{Z}_{m-4}$, and he needs one more step to 
    dominate $v_2$. Hence he needs at least $(m-5)+4=m-1$ steps to win. 
    
    \paragraph{(3.1.2.2.4)} $D_3=u_j$, $j<i$, hence $j\in \{3,4\}$. 
    
    \paragraph{(3.1.2.2.4.a)} If $j=3$, i.e., $D_3=u_3$. Let $S_3=v_1$, which forces
    $D_4=u_2$ by making a line trap. The remaining graph on the right is $\mathfrak{X}_{m-3}$, but with $v_i$ claimed by Dominator already,
    and it is Staller's turn now. By induction, Dominator needs at least $(m-3)-2=m-5$ more steps to win. In total,
    he needs at least $(m-5)+3=m-2$ steps to win. 
    
    \paragraph{(3.1.2.2.4.b)} If $j=4$, i.e., $D_3=u_4$. Let $S_3=u_2$, which forces $D_4= v_1$. Then let $S_4=v_4$,
    which forces $D_5=u_3$ by making a line trap. The remaining graph on the right is $\mathfrak{X}_{m-4}$, 
    but with $v_i$ claimed by Dominator already,
    and it is Staller's turn now. By induction, Dominator needs at least $(m-4)-2=m-6$ more steps to win. In total,
    he needs at least $(m-6)+4=m-2$ steps to win. 
    
     \paragraph{(3.2)} $D_2=v_1$. Let $S_2=u_3$, we carry out case distinctions according to the third 
    move of Dominator in the sequel.
    
    \paragraph{(3.2.1)} $D_3\in\{u_1,u_2\}$. 
    
    \paragraph{(3.2.1.1)} $i$ is even. Let Staller create a sequence of triangle traps $v_3--u_{i-2}$. The remaining 
    graph is the MBD graph $\mathfrak{W}_{m-i}$ on the right, and it is Staller's turn now. By Proposition~\ref{prop:w_m},
    Dominator needs at least $(m-i)-1$ more steps to win. In total, he needs at least $(m-i-1)+(i-1)=m-2$ steps to win. 
    
    \paragraph{(3.2.1.2)} $i$ is odd. Let Staller create a sequence of triangle traps $v_3--u_{i-1}$, which ends with
    $S_{i-1}=u_i$, $D_i=u_{i-1}$. Let $S_i=u_{i+2}$, which creates the MBD graph $\mathfrak{R}_{m-i}$ on the right, and 
    it is Dominator's turn now. By Proposition~\ref{prop:R_m}, $\gamma_{MB}(\mathfrak{R}_{m-i})=m-i$, hence Dominator needs
    at least $m-i$ steps to win. In total, he needs at least $(m-i)+i=m$ moves to win. 
    
    \paragraph{(3.2.2)} $D_3=u_j$, $j\geq 4$. A similar reasoning as in case {\bf (2.1.2.2)}
    shows that $\min\{i,j\}\in\{3,4\}$ must hold. 
    
    \paragraph{(3.2.2.a)} $j<i$. Since $u_3$ is claimed, $j=4$, i.e., $D_3=u_4$. The analysis is analogous to 
    case {\bf (2.4.2.2.a)}. 
    
    \paragraph{(3.2.2.b)} $j=i$. Since $u_3$ is claimed, $i=j=4$, that is, $D_1=v_4$, $D_3=u_4$. Let $S_3=u_2$,
    which forces $D_4=u_1$. Then the remaining graph is the MBD graph $\mathfrak{Z}_{m-4}$ on the right, and it is Staller's 
    turn now. By Proposition~\ref{prop:z_m}, Dominator needs at least $(m-4)-1=m-5$ more steps to win. In total,
    he needs at least $(m-5)+4=m-1$ steps to win. 
    
    \paragraph{(3.2.2.c)} $j>i$, hence $i\in\{3,4\}$. If $i=3$, that is, $D_1=v_3$, then let $S_3=u_2$, which 
    forces $D_4=u_1$. We get the MBD graph $\mathfrak{X}_{m-3}$ on the right, but with $v_j$ claimed already, and it is 
    Staller's turn now. By induction, Dominator needs at least $(m-3)-3=m-6$ more steps to win. 
    In total, he needs at least $(m-6)+4=m-2$ steps to win. If $i=4$, that is, $D_1=v_4$. Let $S_3=u_2$, 
    which forces $D_4=u_1$. Then let $S_4=u_4$, which forces $D_5=v_3$ by creating a line trap. Then consider 
    the MBD graph $\mathfrak{X}_{m-4}$ on the right. Again, Dominator needs in total at least $m-2$ steps to win.

    \paragraph{(3.2.3)} $D_3=v_j$, $j\geq 3$. A similar reasoning as in case {\bf (2.1.2.2)}
    shows that $\min\{i,j\}\in\{3,4\}$ must hold. 
    
    \paragraph{(3.2.3.a)} $j<i$. The reasoning is analogous to cases {\bf (3.2.2.c)}.
    
    \paragraph{(3.2.3.b)} $j>i$. The reasoning is analogous to cases {\bf (3.2.2.c)}.
    
     \paragraph{(3.3)} $D_2=u_j$, where $j\geq 2$, $j>i$ and $i$ is even. Let $S_2=u_2$, the remaining reasoning 
     is analogous to case {\bf (2.7)}.
      
       \paragraph{(3.4)} $D_2=u_j$, where $j\geq 2$, $j>i$ and $i$ is odd. Let $S_2=u_1$, the remaining reasoning 
     is analogous to case {\bf (2.8)}.
       
        \paragraph{(3.5)} $D_2=u_j$, where $j=i\geq 2$. Let $S_2=v_1$, the remaining reasoning 
     is analogous to case {\bf (2.6)}.
        
         \paragraph{(3.6)} $D_2=u_j$, where $2 \leq j<i$. Let $S_2=v_1$, the remaining reasoning 
     is analogous to case {\bf (2.5)}.
         
          \paragraph{(3.7)} $D_2=v_j$, $j\geq 2$, such that $\min\{i,j\}$ is odd. Let $S_2=u_1$, the remaining reasoning 
     is analogous to case {\bf (2.8)}.
          
           \paragraph{(3.8)} $D_2=v_j$, $j\geq 2$, such that $\min\{i,j\}$ is even. Let $S_2=u_2$, the remaining reasoning 
     is analogous to case {\bf (2.7)}.
   
   Since we have gone through all cases, it turns out that Dominator anyway needs at least $m-2$ steps to win the game. 
   By induction, $\gamma_{MB}(\mathfrak{X}_m)\geq m-2$, when $m\geq 6$. So to sum up, 
    \[  \gamma_{MB}(\mathfrak{X}_m)=m-2\; \text{for}\; m\geq 6. \vspace*{-23PT}  \]
\end{itemize}
\end{proof}

\subsection{Maker--Breaker domination number of $P_2\square P_n$}
\begin{proof}[of Theorem~\ref{thm:grid_s}]
To prove the upper bound $\gamma'_{MB}(P_2\square P_n)\leq n$, let Dominator adopt the pairing strategy:
whenever Staller claims $u_i$ ($v_i$), let Dominator claim $v_i$ ($u_i$) in the next round. 
It is not hard to see that Dominator will win within $n$ moves. 

To prove the lower bound $\gamma'_{MB}(P_2\square P_n)\geq n$, we propose the following strategy for Staller.
Let Staller claim $v_2$ as her first move if $n\geq 2$, let Staller claim $u_1$ when $n=1$. 
The remaining game is harder for Dominator than the D-game on $\mathfrak{R}_n$, since 
$\mathfrak{R}_n$ has one less undominated vertex than $P_2\square P_n$ with $v_2$ claimed by Staller. 
By Proposition~\ref{prop:R_m}, we know that Dominator needs at least $n$ more steps to win, and he does not skip any moves.
When $n=1$, it is easy to see that he shall not skip any moves. Hence $\gamma'_{MB}(\mathfrak{R}_n)\geq n$, and Dominator would not win if he 
skipped any moves.
\end{proof}

Next, we would like to address a concrete example of $P_2 \square P_n$, namely when $n=13$. Before which, we need 
two more auxiliary results.
\begin{claim}\label{clm: w'_4}
 Consider the S-game on the MBD graph $\mathfrak{W}_4$. Let Dominator skip his first move, namely 
 $D'_1=\emptyset$. In addition, assume that $S'_1\notin\{u_3,v_3,u_4,v_4\}$. Then Dominator 
 can win within $4$ moves.
\end{claim}
\begin{proof}
 The proof can be done by case distinctions, we do not go into details here.
\end{proof}

\begin{claim}\label{clm:w'_6}
 Consider the S-game on the MBD graph $\mathfrak{W}_6$. Let Dominator skip his first move, namely 
 $D'_1=\emptyset$. In addition, assume that $S'_1= v_2$. Then Dominator can win within $6$ moves.
\end{claim}
\begin{proof}
 The proof can be done by case distinctions, we do not go into details here.
\end{proof}

Now we can prove the following result about $P_2 \square P_{13}$.

\begin{theorem}\label{P2_P13}
$\gamma_{MB}(P_2 \square P_{13}) = 11$. 
\end{theorem}
\begin{proof}
Since $P_2 \square P_{13}$ has one more undominated vertex than $\mathfrak{X}_{13}$, it follows that
\[\gamma_{MB}(P_2 \square P_{13})\geq \gamma_{MB}(\mathfrak{X}_{13}).\]
This is because of the Continuation Principle (Remark~\ref{rem:cont-principle}).
By Theorem~\ref{thm:x_m}, $\gamma_{MB}(\mathfrak{X}_{13})=11$. Hence $\gamma_{MB}(P_2 \square P_{13})\geq 11$. 
Denote by $L=(V_L,E_L)$ the subgraph of $P_2 \square P_{13}$ induced by vertex set 
\[\{u_1,v_1, u_2,v_2,\ldots, u_6,v_6\}.\]
Denote by $R=(V_R,E_R)$ the subgraph of $P_2 \square P_{13}$ induced by vertex set \[\{u_7,v_7, u_8,v_8,\ldots, u_{13},v_{13}\}.\]
For the upper bound, we propose strategy for Dominator. 

Let $D_1=v_7$. It suffices to consider Staller's 
response on $V_L\cup\{u_7\}$, by the symmetric property of $P_2 \square P_{13}$. We carry out case distinctions 
according to the choice of $S_1$ as follows:
\begin{itemize}
 \item { Case 1:}  $S_1=u_7$. In this case, we see an MBD graph $\mathfrak{X}_6$ on both subgraphs $L$ and $R$. Let 
 Dominator carry out the next move on the right MBD graph $\mathfrak{X}_6$. Then, whenever Staller claims some vertex 
 on $R$, let Dominator 
 adopt his strategy on the right MBD graph $\mathfrak{X}_6$; whenever Staller claims some vertex on $L$, 
 let Dominator respond on the left MBD graph $\mathfrak{X}_6$, using the pairing strategy. 
 By Theorem~\ref{thm:x_m}, Dominator needs $4$ steps to dominate all vertices on the right. He needs $6$ steps to dominate all vertices on the left.
 In total, he can win within $6+1+4=11$ steps. 
 
 \item {Case 2:} $S_1=u_5$, then let $D_2=u_9$. Then let Dominator respond on $L$ if Staller claims any vertex of $L$, and let him respond
 on $R$ if Staller claims any vertex of $R$. On the right is the MBD graph $\mathfrak{W}_4$, while on the left it is the situation described in 
 Claim~\ref{clm:w'_6}. By Proposition~\ref{prop:w_m}, Dominator needs $\gamma'_{MB}(\mathfrak{W}_4)= 3$ steps to dominate all vertices of $R$. Claim~\ref{clm:w'_6} tells us that 
 Dominator needs within $6$ steps to dominate all vertices of $L$. Hence, he needs in total within $6+3+2=11$ steps to win the game.
 
 \item {Case 3:} $S_1\in\{u_3,v_3,u_4,v_4,v_5,u_6,v_6\}$, then let $D_2=u_5$. Then we have the MBD graph $\mathfrak{W}_6$ on the right. Let Dominator respond 
 on $V_R\cup\{u_7\}$ whenever Staller claims any vertex in $V_R\cup\{u_7\}$, and let him respond on $V_L$ otherwise. 
 If $S_1\in\{u_6,v_6\}$, on the left would be $\mathfrak{W}_4$. In this case, Dominator needs at most 
 \[\gamma'_{MB}(\mathfrak{W}_4)+\gamma'_{MB}(\mathfrak{W}_6)+2=3+5+2=10\] steps to win.
 Otherwise, $S_1\in\{u_3,v_3,u_4,v_4,v_5\}$, then the graph on the left fulfills the description in Claim~\ref{clm: w'_4}, where he needs within $4$ steps
 to win. In total, he needs at most $4+5+2=11$ steps to win. 
 
\item {Case 4:} $S_1\in\{u_2,v_2\}$, then let $D_2=u_3$. Then let Dominator respond on $V_R\cup\{u_7\}$ whenever Staller claims any vertex of this set.
Let Dominator adopt the pairing strategy on the left, where the pairing sets of vertices are $\{u_1,v_1\}$, $\{v_4,v_5\}$, $\{u_5,u_6\}$. In addition,
he needs at most one more step to dominate $v_2$. 
By Proposition~\ref{prop:w_m}, $\gamma'_{MB}(\mathfrak{W}_6)=5$, hence Dominator needs in total at most $3+1+5+2=11$ steps to win the game. 

\item {Case 5:} $S_1\in\{u_1,v_1\}$, then let $D_2=v_2$. Let Dominator adopt strategy for the S-game on $\mathfrak{W}_6$ on the right, and the pairing strategy 
with pairing sets $\{u_3,u_4\}$, $\{v_4,v_5\}$, $\{u_5,u_6\}$. In addition, he needs at most one more step to dominate $u_1$. In total he needs 
within $3+1+5+2=11$ steps to win the game. 
 
\end{itemize}
To sum up, $\gamma_{MB}(P_2 \square P_{13})\leq 11$. Combining it with the lower bound, we obtain that 
\[ \gamma_{MB}(P_2 \square P_{13})=11.\vspace*{-20PT} \]
\end{proof}

\begin{proof}[of item 3 in Theorem~\ref{thm:grid_d}] 

By Theorem~\ref{P2_P13}, $\gamma_{MB}(P_2\square P_{13})=11$. When $n>13$, consider the graph as the two 
subgraphs $A:=P_2\square P_{13}$
and $B\cong P_2\square P_{n-13}$ connected by edges $\{u_{13},u_{14}\}$, $\{v_{13},v_{14}\}$, up to isomorphism. Let Dominator respond 
on $A$ whenever Staller claims a vertex of $A$, and let Dominator respond on $B$ with the pairing strategy whenever Staller claims a vertex of $B$. 
In this way, we see that he needs within \[11+(n-13)=n-2\] steps in total, in order to win. 
Therefore, \[\gamma_{MB}(P_2\square P_{13})\leq n-2,\; n\geq 13.\]

For the lower bound, by the Continuation Principle (Remark~\ref{rem:cont-principle}) we get 
\[\gamma_{MB}(P_2\square P_{13})\geq \gamma_{MB}(\mathfrak{X}_n),\]
since $P_2\square P_{13}$ has one more undominated vertex than the 
MBD graph $\mathfrak{X}_n$. By Theorem~\ref{thm:x_m}, \[\gamma_{MB}(\mathfrak{X}_n)\geq n-2,\; n\geq 6.\] Hence 
\[\gamma_{MB}(P_2\square P_{13})\geq \gamma_{MB}(\mathfrak{X}_n)\geq n-2,\; n\geq 13.\]
To conclude, 
\[ \gamma_{MB}(P_2\square P_{13})=n-2,\; n\geq 13. \vspace*{-20PT} \]
\end{proof}

\subsection{MBD game on $P_2\square P_{12}$}
In this section, we prove the first two items of Theorem~\ref{thm:grid_d}. We re-state them here as one proposition 
(Proposition~\ref{prop:n<5}) and one theorem (Theorem~\ref{thm:n_5_12}), 
and prove them one by one.
\begin{proposition}\label{prop:n<5}
 If $1\leq n\leq 4$, then $\gamma_{MB}(P_2\square P_n)=n$.
\end{proposition}
\begin{proof}
 For $n=1$ and $n=2$, one can easily do the verification. 
 For $n\in\{3,4\}$, let Dominator adopt the pairing strategy --- we see that $\gamma_{MB}(P_2\square P_n)\leq n$. 
 For the other direction, in the $n=3$ case, because of the symmetry of the graph, it is sufficient to consider two cases,
 namely when $D_1\in\{u_1,u_2\}$.
 \begin{itemize}
  \item Case 1: $D_1=u_1$. Let $S_1 = v_3$. There are still three undominated vertices, it is not hard 
  to see that Dominator needs at least two more moves to win the game.
  \item Case 2: $D_1=u_2$. Let $S_1 = v_2$. One observes that for the remaining two undominated vertices,
  Dominator needs two more moves to win. 
 \end{itemize}
Therefore, $\gamma_{MB}(P_2\square P_3)=3$.

When $n=4$, because of the symmetry of the graph, it is sufficient to consider two cases, namely when $D_1\in\{u_1,u_2\}$.
\begin{itemize}
 \item Case 1: $D_1=u_1$. We see that the remaining graph is the MBD graph~$\mathfrak{W}_3$ on the right. 
 By Proposition~\ref{prop:w_m}, $\gamma'_{MB}(\mathfrak{W}_3)=3$. Hence Dominator needs in total four steps to win.
 \item Case 2: $D_1=u_2$, let $S_1=v_4$. One can verify that Dominator cannot dominate the 
 remaining four vertices within two steps. 
\end{itemize}
Thence, $\gamma_{MB}(P_2\square P_4)=4$.
\end{proof}

\begin{theorem}\label{thm:n_5_12}
 If $5\leq n\leq 12$, then $\gamma_{MB}(P_2\square P_n)=n-1$.
\end{theorem}
In order to prove the above result, we need some preparations.
\begin{lemma}\label{lem:n=5}
 $\gamma_{MB}(P_2\square P_5)=4$.
\end{lemma}
\begin{proof}
 First, we prove $\gamma_{MB}(P_2\square P_5)\geq 4$. 
 It is sufficient to consider three cases because of the symmetry of the graph, namely when $D_1\in \{u_1,u_2,u_3\}$.
 \begin{itemize}
  \item Case 1: $D_1=u_1$. The remaining graph is the MBD graph~$\mathfrak{W}_4$. By Proposition~\ref{prop:w_m}, 
  $\gamma'_{MB}(\mathfrak{W}_4)=3$. Hence Dominator needs in total four steps to win.
  \item Case 2: $D_1=u_2$. Let $S_1=v_2$. We see that the remaining part on the right is the MBD graph~$\mathfrak{X}_3$.
  By Theorem~\ref{thm:x_m}, we have $\gamma_{MB}(\mathfrak{X}_3)=2$. Note that Dominator also needs one step on the left 
  part of the graph, so as to dominate~$v_1$. He needs in total four steps to win.
  \item Case 3: $D_1=u_3$. Let $S_1=v_5$. Then we see that Dominator needs at least two more steps to dominate 
  $v_4,v_5,u_5$. Also he needs at least one more step to dominate $u_1,v_1,v_2$. Hence he needs in total 
  no less than four steps. 
 \end{itemize}
Now we prove the other direction. Let $D_1=u_1$. The remaining graph is the MBD graph~$\mathfrak{W}_4$. 
Hence no matter how Staller plays, Dominator can win within~$3$ steps, since $\gamma'_{MB}(\mathfrak{W}_4)=3$
by Proposition~\ref{prop:w_m}.
\end{proof}

The following several claims are listed here to prepare for the proof of Theorem~\ref{thm:geq_11}.
\begin{claim}\label{claim:v9}
 In the D-game on $P_2\square P_{12}$, if we are in the situation where $D_1=v_6$, $S_1=u_8$, $D_2=u_6$, $S_2=u_9$,
 if $D_3\neq v_9$, then Dominator would not win the game.
\end{claim}
\begin{proof}[of Claim~\ref{claim:v9}]
The analysis of the situation where $D_3\neq v_9$ contains three cases.
\begin{itemize}
\item Case 1: $D_3\notin\{u_7,v_7,v_9,u_{10},v_{10}\}$. Then let $S_3=v_9$ which forces $D_4=u_{10}$. 
Next, let $S_4=v_8$. We see that Dominator is forced to claim three vertices $v_{10}, v_7,u_7$ 
all in the next move, which is not possible. 
\item Case 2: $D_3=u_{10}$ (or $D_3=v_{10}$). Then let $S_3=v_8$ which forces $D_4=u_7$. 
Then let $S_4=v_9$. We see that Dominator must claim two vertices $v_7, v_{10}$ (or $v_7,u_{10}$) in the next round, 
which is obviously not allowed by the rule of this game. 
\item Case 3: $D_3=u_7$ (or $D_3=v_7$). Then let $S_3=v_9$ which forces $D_4=u_{10}$. 
Then let $S_4=v_8$. We see that Dominator must claim two vertices $v_7, v_{10}$ (or $u_7,v_{10}$) in the next round, 
which is obviously not allowed by the rule of this game. \vspace*{-20PT}
\end{itemize}
\end{proof}

\begin{claim}\label{claim2}
 In the D-game on $P_2\square P_{12}$, suppose that we are in the situation where $D_1=v_6$, $S_1=u_8$, $D_2=u_7$, $S_2=v_9$.
 Then if $D_3\notin\{u_9,u_{10},v_{10}\}$, he would not be able to win the game.
\end{claim}
\begin{proof}
 Suppose that $D_3\notin\{u_9,u_{10},v_{10}\}$. We do the analysis in two cases.
 \begin{itemize}
  \item Case 1: $D_3=v_8$. Let Staller make a sequence of triangle traps which means: $S_3=u_{10}\rightarrow D_4=u_9$,
  $S_4=v_{11}\rightarrow D_5=v_{10}$. Next let $S_5=u_{12}$, Dominator could not anymore dominate both~$u_{11}$ and~$v_{12}$.
  \item Case 2: $D_3\notin \{v_8,u_9,u_{10},v_{10}\}$. Let $S_3=u_9$ which forces $D_4=u_{10}$. Next, let $S_4=v_8$.
  Dominator could not anymore dominate both~$v_8$ and~$v_9$. \vspace*{-20PT}
 \end{itemize}
\end{proof}

\begin{claim}\label{claim3}
 In the D-game on $P_2\square P_{12}$, suppose that we are in the situation where $D_1=v_6$, $S_1=u_8$, $D_2=v_7$, $S_2=u_9$.
 Then if $D_3\notin\{v_9,u_{10},v_{10}\}$, he would not win the game.
\end{claim}
\begin{proof}
 Suppose that $D_3\notin\{u_9,u_{10},v_{10}\}$. We do the analysis in two cases.
 \begin{itemize}
  \item Case 1: $D_3=v_8$. Let Staller make a sequence of line traps which means: $S_3=u_{10}\rightarrow D_4=v_9$,
  $S_4=u_{11}\rightarrow D_5=v_{10}$. Next let $S_5=u_{12}$, Dominator would not manage to dominate both~$v_{11}$ and~$v_{12}$ anymore.
  \item Case 2: $D_3\notin \{v_8,u_9,u_{10},v_{10}\}$. Let $S_3=v_8$ which forces $D_4=u_7$. Next, let $S_4=v_9$.
  Dominator would not be able to dominate both~$v_9$ and~$u_9$ anymore. \vspace*{-20PT}
 \end{itemize}
\end{proof}

\begin{claim}\label{claim4}
 In the D-game on $P_2\square P_{12}$, suppose that we are in the situation where $D_1=v_6$, $S_1=u_8$, 
 $D_2\in\{u_i,v_i\mid 10\leq i\leq 12\}$, $S_2=v_8$. Then if $D_3\notin\{u_7,v_7,u_9,v_9\}$, Dominator would not be able to win the game.
\end{claim}
\begin{proof}
 Suppose that $D_3\notin\{u_7,v_7,u_9,v_9\}$. We do the analysis in two cases.
 \begin{itemize}
  \item Case 1: $D_3=u_6$. Then \[S_3=v_9\rightarrow D_4=v_7,\; S_4=u_9\rightarrow D_5=u_7.\] Now at least one 
  of~$u_{10}$ and~$v_{10}$ is still free. When~$u_{10}$ (or~$v_{10}$) is free. Then let $S_5=u_{10}$ 
  (or $S_5=v_{10}$), we see that~$u_9$ (or~$v_9$) would be isolated. 
  \item Case 2: $D_3\neq u_6$. Then $S_3=u_7\rightarrow D_4=u_9$. Next let $S_4=v_7$. Then at least one of~$u_7$
  and~$v_8$ would not anymore be dominated in the game. \vspace*{-20PT}
 \end{itemize}
\end{proof}

\begin{claim}\label{claim5}
 In the D-game on $P_2\square P_{12}$, suppose that we are in the situation where $D_1=v_3$, $S_1=u_5$, 
 $D_2=u_6$, $S_2=v_8$, $D_3=v_6$, $S_3=u_8$. Then if $D_4\notin\{u_9,v_9\}$, Dominator would not win the game.
\end{claim}
\begin{proof}
 The proof is left as an exercise. 
\end{proof}

\begin{figure}
\centering
\includegraphics[width=0.742\linewidth]{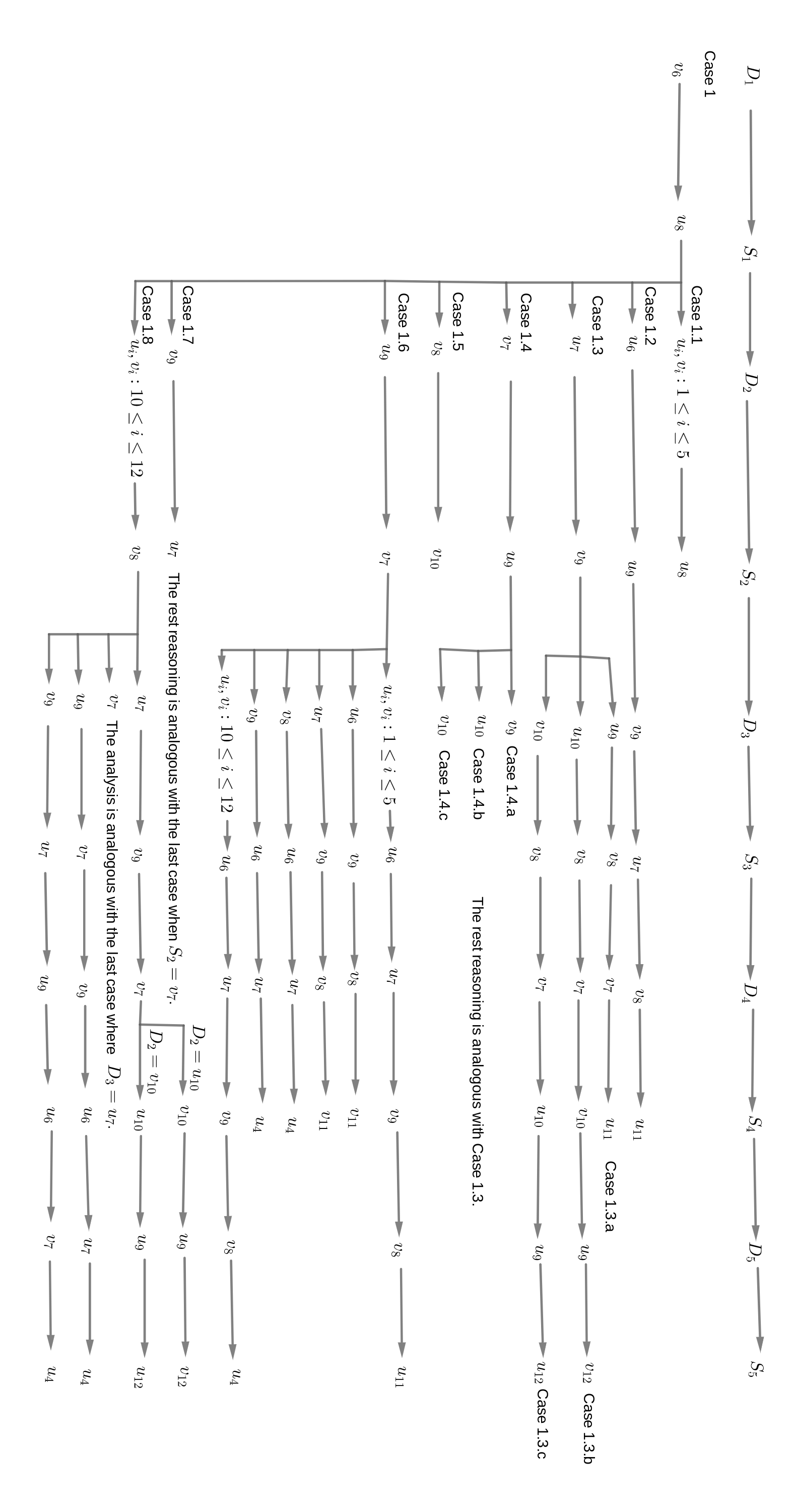}
\caption{This is a strategy tree showing the strategies for Staller under the situation when $D_1=v_6$, in the Maker--Breaker
domination game on $P_2\square P_{12}$. One can verify that with the provided strategy for Staller, Dominator always
needs at least~$11$ steps to win.}
\label{fig:strategy_tree_2}
\end{figure}

The following result, as the last preparation needed for the proof of Theorem~\ref{thm:n_5_12}, has a rather lengthy proof. 
In order to assist our readers to attain a better understanding,
we provide the strategy tree ( see Figure~\ref{fig:strategy_tree_2} ) for the analysis of Case~1 of the proof of Theorem~\ref{thm:geq_11}.

 \begin{theorem}\label{thm:geq_11}
$\gamma_{MB}(P_2 \square P_{12}) \geq 11$.
\end{theorem}
\begin{proof}
 Because of the symmetry of the graph, it is sufficient to consider the cases when 
\[D_1 \in \{v_1, v_2, v_3, v_4, v_5, v_6\}.\] 
\begin{itemize}
\item Case 1: $D_1 = v_6$. Let $S_1=u_8$. 

\subitem Case 1.1: $D_2\in\{u_i,v_i\mid 1\leq i\leq 5\}$. In this case, let $S_2=u_8$. We see that 
on the right part of the graph is the MBD graph~$\mathfrak{R}_6$. We let Staller plays on this part until 
Dominator dominates all vertices on this~$\mathfrak{R}_6$. By Proposition~\ref{prop:R_m}, we know 
that he needs~$6$ steps to accomplish this task and he shall not skip any moves in the middle of this 
process. After Dominator's last move on this~$\mathfrak{R}_6$, let Staller continue her response towards~$D_1$ on 
the left subgraph which is the MBD graph~$\mathfrak{X}_5$. To view it in another way, it is an MBD game
on~$\mathfrak{X}_5$. By Theorem~\ref{thm:x_m}, $\gamma_{MB}(M_5)=4$. So, in total, Dominator needs $1+6+5=12$ moves.

\subitem Case 1.2: $D_2=u_6$. Then let $S_2=u_9$. In this case, Dominator is forced to claim~$v_9$ ---
for the reasoning see Claim~\ref{claim:v9}. Then let $S_3=u_7$ which forces $D_4=v_8$. 
Next, we let $S_4=u_{11}$. On the left is the MBD graph~$\mathfrak{R}_3$, on the right is the MBD graph~$\mathfrak{Z}_5$.
We let Staller play on~$\mathfrak{R}_3$ until the game finishes on this part, and then let her starts playing
on~$\mathfrak{Z}_5$. By Proposition~\ref{prop:R_m} and~\ref{prop:z_m}, $\gamma_{MB}(\mathfrak{R}_3)=3$ and Dominator would not
skip any moves, $\gamma'_{MB}(\mathfrak{Z}_5)=4$. Hence, Dominator needs in total $4+3+4=11$ moves.

\subitem Case 1.3: $D_2=u_7$. Let $S_2=v_9$. Now by Claim~\ref{claim2}, we know that $D_3\in\{u_9,u_{10},v_{10}\}$.

\subsubitem Case 1.3.a: $D_3=u_9$. Let $S_3=v_8$ which forces $D_4=v_7$. Next, let $S_4=u_{11}$.
We see~$\mathfrak{W}_5$ on the left and~$\mathfrak{R}_3$ on the right. Let Staller plays on~$\mathfrak{R}_3$ first. Since 
$\gamma'_{MB}(\mathfrak{W}_5)=4$ and $\gamma_{MB}(\mathfrak{R}_3)=3$, Dominator needs $4+3+4=11$ moves to complete the game. 

\subsubitem Case 1.3.b: $D_3=u_{10}$. Let $S_3=v_8$ which forces $D_4=v_7$. Then, let $S_4=v_{10}$
which forces $D_5=u_9$. Next, let $S_5=v_{12}$. We see that the remaining part on the right is the 
MBD graph~$\mathfrak{R}_2$, while that on the left is the MBD graph~$\mathfrak{W}_5$. 
We let Staller first play on~$\mathfrak{R}_2$ until the game finishes on this part, and then let her starts playing
on~$\mathfrak{W}_5$. By Proposition~\ref{prop:R_m} and~\ref{prop:w_m},
we know that $\gamma_{MB}(\mathfrak{R}_2)=2$ and Dominator shall not skip any moves during the process, $\gamma'_{MB}(\mathfrak{W}_5)=4$. 
Hence, Dominator needs in total $5+2+4=11$ moves to accomplish the victory. 

\subsubitem Case 1.3.c: $D_3=v_{10}$. Let $S_3=v_8$ which forces $D_4=v_7$. Then, let $S_4=u_{10}$
which forces $D_5=u_9$. Next, let $S_5=u_{12}$. We see that the remaining part on the right is the 
MBD graph~$\mathfrak{R}_2$, while that on the left is the MBD graph~$\mathfrak{W}_5$. 
We let Staller first play on~$\mathfrak{R}_2$ until the game finishes on this part, and then let her starts playing
on~$\mathfrak{W}_5$. By Proposition~\ref{prop:R_m} and~\ref{prop:w_m},
we know that $\gamma_{MB}(\mathfrak{R}_2)=2$ and Dominator shall not skip any moves during the process, $\gamma'_{MB}(\mathfrak{W}_5)=4$. 
Hence, Dominator needs in total $5+2+4=11$ moves to accomplish the victory. 

\subitem Case 1.4: $D_2=v_7$. Let $S_2=u_9$. This case is analogous to Case 1.3. For the completeness
of the contents, we still illustrate the detailed argument.
By Claim~\ref{claim3} we know that $D_3\in\{v_9,u_{10},v_{10}\}$.

\subsubitem Case 1.4.a: $D_3=v_9$. Let $S_3=v_8$ which forces $D_4=u_7$. Then let $S_4=u_{11}$ which then creates
an~$\mathfrak{R}_3$ on the right and let Staller play on this~$\mathfrak{R}_3$ first and then play on the leftmost~$\mathfrak{W}_5$. Dominator needs 
$4+3+4=11$ moves to win.

\subsubitem Case 1.4.b: $D_3=u_{10}$. Let $S_3=v_8$ which forces $D_4=u_7$. Then let $S_4=v_{10}$ which forces 
$D_5=v_9$. Then let $S_5=v_{12}$. Analogous to the cases before, it is~$\mathfrak{R}_2$ on the right and~$\mathfrak{W}_5$ on the left. 
Let Staller play on~$\mathfrak{R}_2$ first and then on~$\mathfrak{W}_5$. Dominator needs $5+2+4=11$ moves to win.

\subsubitem Case 1.4.c: $D_3=v_{10}$. Adopt the following strategy, where~$\rightarrow$ refers 
to ``forces'' (and this convention applies also to the later context of this paper):
\[S_3=u_{10}\rightarrow D_4=v_9,\; S_4=u_7\rightarrow D_5=v_8.\] 
Let $S_5= u_{12}$ and let Staller play on the rightmost~$\mathfrak{R}_2$ until the game is locally finished then play on
the rightmost~$\mathfrak{W}_5$. Dominator needs $5+2+4=11$ moves to win.

\subitem Case 1.5: $D_2=v_8$. Let $S_2=u_{10}$. It is an~$\mathfrak{R}_4$ on the rightmost. By Proposition~\ref{prop:R_m},
$\gamma_{MB}(\mathfrak{R}_4)=4$ and Dominator shall not skip any moves in the process. After the game is finished on this~$\mathfrak{R}_4$,
let Staller then claim~$u_6$. We get an~$\mathfrak{X}_5$ on the leftmost. Since Dominator needs one more step to dominate~$u_7$
and $\gamma_{MB}(\mathfrak{X}_5)=4$ (by Theorem~\ref{thm:x_m}), Dominator needs in total $2+4+4+1=11$ steps to win. 

\subitem Case 1.6: $D_2=u_9$. Let~$S_2=v_7$. Consider the following cases for Dominator's third move.

\subsubitem Case 1.6.a: $D_3\in \{u_i, v_i\mid 1\leq i\leq 5\}$. Then 
\[S_3=u_6\rightarrow D_4=u_7,\; S_4=v_9\rightarrow D_5=v_8.\] Then let $S_5= v_{11}$ which creates an~$\mathfrak{R}_3$ on the rightmost. 
Let Staller plays on this~$\mathfrak{R}_3$ first and then starts playing on the left~$\mathfrak{X}_5$ --- when forgetting
the move~$D_3$, it is an~$\mathfrak{X}_5$. Note that~$D_3$ may not be the optimal move, hence Dominator 
needs at least $\gamma_{MB}(\mathfrak{X}_5)=4$ (by Theorem~\ref{thm:x_m}) to finish the local game on the 
leftmost. And on the leftmost~$\mathfrak{R}_3$, he needs~$3$ moves and shall not skip any moves by Proposition~\ref{prop:R_m}.
Hence he needs in total $4+3+4=11$ moves to win the whole game.

\subsubitem Case 1.6.b: $D_3=u_6$. Let $S_3=v_9\rightarrow D_4=v_8$. Then let $S_4=v_{11}$ which 
then creates an~$\mathfrak{R}_3$ on the rightmost, while on the leftmost is an~$\mathfrak{Z}_5$. Similarly as always, 
let us just omit this repetitive narrative from now on. Dominator needs in total $4+3+4=11$ steps to win 
( since $\gamma'_{MB}(\mathfrak{Z}_5)=4$ ).

\subsubitem Case 1.6.c: $D_3=u_7$. Let $S_3=v_9\rightarrow D_4=v_8$. Then $S_4=v_{11}$. Rightmost:~$\mathfrak{R}_3$;
leftmost:~$\mathfrak{W}_5$. Steps to win: $4+3+4=11$.

\subsubitem Case 1.6.d: $D_3=v_8$. Then $S_3=u_6\rightarrow D_4=u_7$, then $S_4=u_4$. Leftmost:~$\mathfrak{R}_5$;
rightmost:~$\mathfrak{W}_3$. Steps to win: $4+5+3=12$.

\subsubitem Case 1.6.e: $D_3=v_9$. Then $S_3=u_6\rightarrow D_4=u_7$, then $S_4=u_4$. Leftmost:~$\mathfrak{R}_5$;
rightmost:~$\mathfrak{Z}_3$. Steps to win: $4+5+2=11$.

\subsubitem Case 1.6.f: $D_3\in\{u_i,v_i\mid 10\leq i\leq 12\}$. Then 
\[S_3=u_6\rightarrow D_4=u_7,\; S_4=v_9\rightarrow D_5=v_8.\] Then let $S_5=u_4$ which creates an~$\mathfrak{R}_5$ on the leftmost. 
Let Staller plays on this~$\mathfrak{R}_5$ first and then starts playing on the rightmost~$\mathfrak{X}_3$ --- when forgetting
the move~$D_3$, it is an~$\mathfrak{X}_3$. Note that~$D_3$ may not be the optimal move, hence Dominator 
needs at least $\gamma_{MB}(\mathfrak{X}_3)=2$ (by Theorem~\ref{thm:x_m}) to finish the local game on the 
leftmost. And on the leftmost~$\mathfrak{R}_5$, he needs~$5$ moves and shall not skip any moves by Proposition~\ref{prop:R_m}.
Hence he needs in total $4+5+2=11$ moves to win the whole game.

\subitem Case 1.7: $D_2=v_9$. Let $S_2=u_7$. The proof is analogous to Case 1.6. We omit it here.

\subitem Case 1.8: $D_2\in\{u_i,v_i\mid 10\leq i\leq 12\}$. Let $S_2=v_8$. Now, by Claim~\ref{claim4}, 
\[D_3\in\{u_7,v_7,u_9,v_9\}.\] 

\subsubitem Case 1.8.a: $D_3=u_7$. Then $S_3=v_9\rightarrow D_4=v_7$. Now we look back on the choice
of~$D_2$: If $D_2\notin\{u_{10},v_{10}\}$, then he cannot win --- because otherwise we let $S_4=u_9$
which forces Dominator in the next step to claim both~$u_{10}$ and~$v_{10}$ (in order to not 
isolate~$u_9$ and~$v_9$) which is not allowed by the rules of this game. Therefore we only need to
consider the situation where Dominator actually has an opportunity to win.
\begin{itemize}
\item Situation 1: $D_2=u_{10}$. Then $S_4=v_{10}\rightarrow D_5=u_9$. Next let $S_5=v_{12}$. 
Now it is back to our familiar pattern. Leftmost:~$\mathfrak{R}_2$; rightmost:~$\mathfrak{W}_5$. Since $\gamma_{MB}(\mathfrak{R}_2)=2$ 
and Dominator shall not skip any moves in the game, $\gamma'_{MB}(\mathfrak{W}_5)=4$, he needs $5+2+4=11$ moves to win. 
\item Situation 2: $D_2=v_{10}$. Then $S_4=u_{10}\rightarrow D_5=u_9$. Next let $S_5=u_{12}$. 
Now it is back to our familiar pattern. Leftmost:~$\mathfrak{R}_2$; rightmost:~$\mathfrak{W}_5$. Since $\gamma_{MB}(\mathfrak{R}_2)=2$ 
and Dominator shall not skip any moves in the game, $\gamma'_{MB}(\mathfrak{W}_5)=4$, he needs $5+2+4=11$ moves to win. 
\end{itemize}

\subsubitem Case 1.8.b: $D_3=v_7$. The reasoning is analogous to Case 1.8.a --- we omit it here. 

\subsubitem Case 1.8.c: $D_3=u_9$. Then 
\[S_3=v_7\rightarrow D_4=v_9,\; S_4=u_6\rightarrow D_5=u_7.\]
Then let $S_5=u_4$ which creates an~$\mathfrak{R}_5$ on the leftmost, and on the rightmost is an~$\mathfrak{Z}_3$ if 
we forget the move of~$D_2$. Since~$D_2$ may not be the optimal choice, Dominator needs at
least $\gamma'_{MB}(\mathfrak{Z}_3)=2$ steps to finish the game on the rightmost subgraph~$\mathfrak{Z}_3$. In total, at least
$4+5+2=11$ moves are needed to win.

\subsubitem Case 1.8.d: $D_3=v_9$. Then 
\[S_3=u_7\rightarrow D_4=u_9,\; S_4=u_6\rightarrow D_5=v_7.\] 
Next let $S_5=u_4$. The rest of the reasoning is the same as in Case 1.8.c.

\item Case 2: $D_1=v_5$. In this case, the analysis is analogous to in Case 1 --- the only change is that
when Dominator (or Staller) claims vertex~$u_i$ (or~$v_i$) in Case 1, they instead now in this case claims
vertex~$u_{i-1}$ (or~$v_{i-1}$), where $i\geq 1$. For example, we already see that then we have $S_1=u_7$ (it was~$u_8$ in Case~1).
Note that this is feasible since all the vertices claimed by the two players in Case 1 has subscript bigger than~$3$.
We have verified all the cases, and we omit the detailed reasoning (or rather repetition) here.

\item Case 3: $D_1=v_4$. In this case, the analysis is analogous to in Case 1 --- the only change is that
when Dominator (or Staller) claims vertex~$u_i$ (or~$v_i$) in Case 1, they instead now in this case claims
vertex~$u_{i-2}$ (or~$v_{i-2}$), where $i\geq 1$. For example, we already see that then we have $S_1=u_6$ 
(it was~$u_8$ in Case~1). Note that this is feasible since all the vertices claimed by the two players 
in Case~1 has subscript bigger than~$3$.
We have verified all the cases, and we omit the detailed reasoning here.

\item Case 4: $D_1=v_3$. Let $S_1=u_5$. In this case, most of the the analysis is analogous to in Case 1 --- the only change is that
when Dominator (or Staller) claims vertex~$u_i$ (or~$v_i$) in Case 1, they instead now in this case claims
vertex~$u_{i-3}$ (or~$v_{i-3}$), where $i\geq 1$. Note that this is feasible since all the vertices claimed by the two players 
in Case~1 has subscript bigger than~$3$. But there are two cases we need to consider separately, namely 
when $D_2\in\{u_6,v_6\}$ --- because, if in these two cases we let Staller adopt the analogous strategy as in Case~1,
Dominator would be able to win within~$10$ moves. We have verified all other cases. Now, let us talk about these
two cases.

\subitem Case 4.1: $D_2=u_6$. Let $S_2=v_8$. Now, Staller's strategy depends on Dominator's third move.

\subsubitem Case 4.1.1: $D_3\in\{u_i,v_i\mid 1\leq i\leq 5\}$. Let $S_3=v_6$ 
which then creates an~$\mathfrak{R}_6$ on the rightmost. Let Staller play on this~$\mathfrak{R}_6$ now.
By Proposition~\ref{prop:R_m}, Dominator needs~$6$ steps to win and he 
shall not skip any moves during the process. Next, Staller's strategy depends on 
what was~$D_3$:
\begin{itemize}
\item Case 4.1.1.a: $D_3\in\{u_1,v_1,u_2,v_2\}$. Let $S_9= v_4$ which forces $D_{10}=v_5$. Then $S_{10}=u_3\rightarrow D_{11}=u_4$. 
Dominator needs at least $4+6+1=11$ steps to win.
\item Case 4.1.1.b: $D_3=u_3$. Then $S_9=v_4\rightarrow D_{10}=v_5$. One more step is needed on the leftmost to dominate~$u_1$ and~$v_1$.
He needs in total $4+6+1=11$ moves to win.
\item Case 4.1.1.c: $D_3=u_4$. Then $S_9=v_5\rightarrow D_{10}=v_4$. On the left is an~$\mathfrak{W}_2$. Steps to win: $4+6+2=12$.
\item Case 4.1.1.d: $D_3=v_4$. Then $S_9=u_1$ which forces Dominator to spend two more steps to finish the game. 
Steps to win: $3+6+2=11$. 
\item Case 4.1.1.e: $D_3=v_5$. Then it is not hard to see that he needs in the best case two steps to win, so he needs at 
least $3+6+2=11$ steps to finish the game. 
\end{itemize}
\subsubitem Case 4.1.2: $D_3=v_9$. Let $S_3=u_8$. By Claim~\ref{claim5}, $D_4\in\{u_9,v_9\}$. We leave the rest 
of the argument for readers as an exercise. 

From now on, we only list the cases, since we believe after all the case distinctions above, our readers already 
got the taste of how this can be done. We have verified all cases, but we omit the detailed long-winded here. 

\subsubitem Case 4.1.3: $D_3=u_7$. Let $S_3=v_9$. Then it is forced that $D_4\in\{u_9,u_{10},v_{10}\}$.

\subsubitem Case 4.1.4: $D_3=v_7$. Let $S_3=u_9$. Then it is forced that $D_4\in\{u_9,u_{10},v_{10}\}$.

\subsubitem Case 4.1.5: $D_3=u_8$. Let $S_3=v_{10}$. 

\subsubitem Case 4.1.6: $D_3=u_9$. Let $S_3=v_7$. 

\subsubitem Case 4.1.7: $D_3=v_9$. Then $S_3=u_7$.

\subsubitem Case 4.1.8: $D_3\in\{u_{10},u_{11},u_{12},v_{10},v_{11},v_{12}\}$. Let $S_3=u_8$. Then it has to be that
$D_4\in\{u_9,v_9\}$.

\subitem Case 4.2: $D_2=v_6$. The reasoning is analogous to Case 4.1.

\item Case 5: $D_1=v_2$. Then $S_1=u_4$. In this case, the analysis is analogous to in Case 4 --- the only change is that
when Dominator (or Staller) claims vertex~$u_i$ (or~$v_i$) in Case 4, they instead now in this case claims
vertex~$u_{i-1}$ (or~$v_{i-1}$), where $i\geq 2$. 
Note that in the case when the subscript of the claimed vertex is~$1$, we keep it unchanged.
We have verified all the cases, and we omit the detailed reasoning (or rather repetition) here.

\item Case 6: $D_1=v_1$. The rest is the MBD graph~$\mathfrak{W}_{11}$. By Proposition~\ref{prop:w_m}, $\gamma'_{MB}(\mathfrak{W}_{11})=10$.
So Dominator needs in total $1+10=11$ steps to win. \vspace*{-20PT}
\end{itemize}
\end{proof}

\begin{proof}[of Theorem~\ref{thm:n_5_12}]
 Let $n\geq 5$, we consider the D-game on $P_2\square P_n$. Consider the graph as two subgraphs joined by 
 two edges: $P_2\square P_5$ on the leftmost and the rightmost is isomorphic to $P_2\square P_{n-5}$; obviously
 they are joined by two edges $\{u_5,u_6\}$, and $\{v_5,v_6\}$. Let Dominator start the game on the leftmost 
 $P_2\square P_5$. Then, whenever Staller claims a vertex on this part, he responds using his strategy 
 for $P_2\square P_5$. Otherwise, let him adopt the pairing strategy to respond to Staller's move. 
 By Lemma~\ref{lem:n=5}, $\gamma_{MB}(P_2\square P_5)=4$. Dominator can win the game within $4+(n-5)=n-1$ steps.
 Therefore, \[\gamma_{MB}(P_2\square P_n)\leq n-1\; \text{when}\; n\geq 5.\] 
 
 Let $5<m<12$. Suppose that $\gamma_{MB}(P_2\square P_m)<m-1$. Then consider the MBD game on $P_2\square P_{12}$.
 We consider the graph as two subgraphs joined by two edges: $P_2\square P_m$ on the leftmost and the rightmost is
 isomorphic to $P_2\square P_{12-m}$. Let Dominator starts playing on the leftmost subgraph. 
 Then, whenever Staller claims a vertex on this part, he responds using his strategy 
 for $P_2\square P_m$. Otherwise, let him adopt the pairing strategy to respond to Staller's move.
 We see that Dominator will win the game with less than $(m-1)+(12-m)=11$ steps, which contradicts 
 the result of Theorem~\ref{thm:geq_11}. 
 Therefore, 
 \[\gamma_{MB}(P_2\square P_n)\geq n-1,\; \text{for}\; 5\leq n\leq 12.\] 
 
 In the end, we get that $\gamma_{MB}(P_2\square P_n)= n-1$, for $5\leq n\leq 12$. 
\end{proof}

\subsection{Union of grids}\label{subsec:union}

In this subsection, we prove the two results on the disjoint union of $P_2\square P_n$s, using Theorem~\ref{thm:grid_s} 
and Theorem~\ref{thm:grid_d}.

\begin{proof}[of Theorem~\ref{thm:union_d}] 
 By the pairing strategy where the pairing sets are $\{u_i,v_i\}$s in each copy, Dominator can win within 
 $k\cdot n$ steps. However, by Theorem~\ref{thm:grid_s}, we know that no matter in which copy Staller claims a vertex,
 Dominator has to respond on the same copy, otherwise he would lose the game. Hence he needs 
 at least $k\cdot \gamma'_{MB}(P_2\square P_n)=k\cdot n$ steps to win. So we get 
 \[
 \gamma_{MB}'(\dot\cup_{i=1}^k(P_2\square P_n)_i)=k\cdot n,\; \text{for any}\; k,n\geq 1. \vspace*{-20PT} \]
\end{proof}

\begin{proof}[of Theorem~\ref{thm:union_s}] 
First we prove the upper bounds. Let Dominator adopt the following strategy: let him respond on the same copy where Staller played her 
move in the last round. 
 By Theorem~\ref{thm:grid_s}, we know that Dominator does not skip any moves in a copy where Staller starts the game --- namely, he has to respond 
 on the same copy whenever Staller plays a move on that copy --- and needs~$n$ steps to 
 dominate all vertices of this copy. For the copy where he starts the game, the situation is just a D-game on $P_2\square P_n$. 
 Hence we get that 
 \[\gamma_{MB}(\dot\cup_{i=1}^k(P_2\square P_n)_i)\leq (k-1)\cdot n +\gamma_{MB}(P_2\square P_n).\]
 
 Now we prove the lower bounds. After Dominator's first move, let Staller continue the game on any other copy than the one where Dominator played 
 his first move, and let Staller continue play on this copy using the strategy for $P_2\square P_n$ until the game is finished on this copy.
 By Theorem~\ref{thm:grid_s}, we know that Dominator does not skip any moves during this process and needs~$n$ steps to 
 dominate all vertices of this copy.
 Then, let Staller repeat this strategy on all the remaining copies except for the one where Dominator starts the game.
 After this process, let Staller continue by responding on that remaining copy, then it should be considered 
 as a D-game on $P_2\square P_n$. 
 Hence we get that 
 \[\gamma_{MB}(\dot\cup_{i=1}^k(P_2\square P_n)_i)\geq (k-1)\cdot n +\gamma_{MB}(P_2\square P_n).\]
 
 Combining the above results, we get that 
 \[\gamma_{MB}(\dot\cup_{i=1}^k(P_2\square P_n)_i)= (k-1)\cdot n +\gamma_{MB}(P_2\square P_n).\]
 By Theorem~\ref{thm:grid_d}, we get the results in the statement of this theorem.
\end{proof}

\section{Concluding remarks}

We proved that Dominator needs exactly~$n$ moves to win in the S-game on $P_2 \square P_n$ for every $n\geq 1$, 
while in the D-game he needs exactly~$n$, $n-1$, $n-2$ moves for $1\leq n\leq 4$, $5\leq n\leq 12$, and $n\geq 13$, respectively.
We showed that for~$k$ copies of $P_2\square P_n$, Dominator needs exactly $k\cdot n$ 
steps to win in the S-game, and $k\cdot n$, $k\cdot n-1$, $k\cdot n-2$ steps to win in the D-game
for $1\leq n$, $5\leq n\leq 12$, $n\leq 13$, respectively. 
The exact result for general Cartesian $P_m\square P_n$ does not seem easy, so it would be interesting to consider 
the situation for $P_3\square P_n$, as a starting point. 
\newline 
\newline
\newline
\acknowledgements

The authors would like to thank Professor Mirjana Mikala\v{c}ki for the helpful comments which contributed to the improvement of the paper. 
The second author thanks Professor Josef Schicho for the general scientific research guidance during this work. The second author thanks 
Dongsheng Wu\footnote{At the time of the publication of this paper, Dongsheng Wu was affiliated with Yanqi Lake Beijing Institute of Mathematical 
Sciences and Applications (BIMSA), and with Tsinghua University, Beijing, China. } for several times' inspiring discussion,
especially on the Maker--Breaker games on the graph $P_2\square P_{13}$.

\nocite{*}
\bibliographystyle{plainnat}
\bibliography{sample-dmtcs}

\end{document}